\documentclass[twoside,reqno]{amsart}

\usepackage[normalem]{ulem}
\usepackage[square,comma]{natbib}

\usepackage[all]{xy}
\xyoption{tips}
\SelectTips {xy} {12}

\usepackage{amsmath,amsthm,amscd,upref,indentfirst}
\usepackage{amsfonts,mathrsfs}
\usepackage{amssymb,amsbsy,bm}

\usepackage{graphicx}
\usepackage[UglyObsolete,small,nohug,heads=LaTeX]{diagrams}
\diagramstyle[labelstyle=\scriptstyle]

\theoremstyle{plain}
\newtheorem{theorem}{Theorem}
\newtheorem{assertion}[theorem]{Assertion}

\theoremstyle{definition}
\newtheorem{definition}[theorem]{Definition}

\theoremstyle{remark}
\newtheorem{remark}[theorem]{Remark}

\newtheorem{example}[theorem]{Example}

\numberwithin{equation}{section}
\numberwithin{theorem}{section}

\renewcommand{\mathfrak}{\textsc}
\renewcommand{\mathcal}{\mathscr}
\renewcommand{\mathbf}{\bm}

\normalfont\upshape  

\title{{A ``\lowercase{q}-deformed'' generalization of the Hossz\'u-Gluskin theorem}
}
\author{{Steven Duplij}}
\address{Mathematisches Institut,
Universit\"at M\"unster,
Einsteinstr. 62,
D-48149 M\"unster,
Germany} 
\email{duplijs@math.uni-muenster.de, sduplij@gmail.com}
\urladdr{http://wwwmath.uni-muenster.de/u/duplij}
\subjclass[2010]{08B05, 17A42, 20N15}

\begin{document}

\thispagestyle{empty}

\begin{abstract}
In this paper a new form of the Hossz\'u-Gluskin theorem is presented
 in terms of polyadic powers and using the language of diagrams. 
It is shown that the Hossz\'u-Gluskin chain formula is 
not unique and can be generalized (``deformed'') using a parameter 
q which takes special integer values. 
A version of the ``q-deformed'' analog of the 
Hossz\'u-Gluskin theorem in the form of an invariance is formulated, 
and some examples are considered. 
The ``q-deformed'' homomorphism theorem is also given.

\end{abstract}

\maketitle

\tableofcontents

\section{Introduction}

Since the early days of \textquotedblleft polyadic
history\textquotedblright\ \cite{kas,pru,dor3}, the interconnection between
polyadic systems and binary ones has been one of the main areas of interest
\cite{pos,thu54}. Early constructions were confined to building some special
polyadic (mostly ternary \cite{cer,leh32}) operations on elements of binary
groups \cite{mil35,leh36,chu}. A very special form of $n$-ary multiplication
in terms of binary multiplication and a special mapping as a chain formula was
found in \cite{hos} and \cite{glu64,glu1}. The theorem that any
$n$-ary multiplication can be presented in this form is called the
Hossz\'{u}-Gluskin theorem (for review see \cite{dud/gla,galmak1}). A concise
and clear proof of the Hossz\'{u}-Gluskin chain formula was presented in
\cite{sok}.

In this paper we give a new form of the Hossz\'{u}-Gluskin theorem in
terms of polyadic powers. Then we show that the Hossz\'{u}-Gluskin chain
formula is not unique and can be generalized (\textquotedblleft
deformed\textquotedblright) using a parameter $q$ which takes special integer
values. We present the \textquotedblleft%
$q$-deformed\textquotedblright\ analog of the Hossz\'{u}-Gluskin theorem in the form of an invariance and
consider some examples. The \textquotedblleft$q$-deformed\textquotedblright%
\ homomorphism theorem is also given.

\section{Preliminaries}

We will use the concise notations from our previous review paper
\cite{dup2012}, while here we repeat some necessary definitions using the
language of diagrams. For a non-empty set $G$, we denote its elements by
lower-case Latin letters $g_{i}\in G$ and the $n$\textit{-tuple} (or
\textit{polyad}) $g_{1},\ldots,g_{n}$ will be written by $\left(  g_{1}%
,\ldots,g_{n}\right)  $ or using one bold letter with index $\mathbf{g}%
^{\left(  n\right)  }$, and an $n$-tuple with equal elements by $g^{n}$. In
case the number of elements in the $n$-tuple is clear from the context or is
not important, we denote it in one bold letter $\mathbf{g}$ without indices.
We omit $g\in G$, if it is obvious from the context.

The Cartesian product $\overset{n}{\overbrace{G\times\ldots\times G}%
}=G^{\times n}$ consists of all $n$-tuples $\left(  g_{1},\ldots,g_{n}\right)
$, such that $g_{i}\in G$, $i=1,\ldots,n$. The $i$\textit{-projection} of the
Cartesian product $G^{n}$ on its $i$-th \textquotedblleft
axis\textquotedblright\ is the map $\mathsf{Pr}_{i}^{\left(  n\right)
}:G^{\times n}\rightarrow G$ such that $\left(  g_{1},\ldots g_{i}%
,\ldots,g_{n}\right)  \longmapsto g_{i}$. The $i$\textit{-diagonal}
$\mathsf{Diag}_{n}:G\rightarrow G^{\times n}$ sends one element to the equal
element $n$-tuple $g\longmapsto\left(  g^{n}\right)  $. The one-point set
$\left\{  \bullet\right\}  $ is treated as a unit for the Cartesian product,
since there are bijections between $G$ and $G\times\left\{
\bullet\right\}  ^{\times n}$, where $G$ can be on any place. In diagrams, if
the place is unimportant, we denote such bijections by $\epsilon$. On the
Cartesian product $G^{\times n}$ one can define a polyadic ($n$-ary or $n$-adic,
if it is necessary to specify $n$, its arity or rank) operation $\mu
_{n}:G^{\times n}\rightarrow G$. For operations we use small Greek letters and
place arguments in square brackets $\mu_{n}\left[  \mathbf{g}\right]  $. The
operations with $n=1,2,3$ are called \textit{unary, binary and ternary}. The
case $n=0$ is special and corresponds to fixing a distinguished element of
$G$, a \textquotedblleft constant\textquotedblright\ $c\in G$, and it is
called a \textit{0-ary operation} $\mu_{0}^{\left(  c\right)  }$, which maps
the one-point set $\left\{  \bullet\right\}  $ to $G$, such that $\mu
_{0}^{\left(  c\right)  }:\left\{  \bullet\right\}  \rightarrow G$, and
(formally) has the value $\mu_{0}^{\left(  c\right)  }\left[  \left\{
\bullet\right\}  \right]  =c\in G$. The composition of $n$-ary and $m$-ary
operations $\mu_{n}\circ\mu_{m}$ gives a $\left(  n+m-1\right)  $-ary operation
by the iteration $\mu_{n+m-1}\left[  \mathbf{g},\mathbf{h}\right]  =\mu
_{n}\left[  \mathbf{g},\mu_{m}\left[  \mathbf{h}\right]  \right]  $. If we
compose $\mu_{n}$ with the 0-ary operation $\mu_{0}^{\left(  c\right)  }$,
then we obtain the arity \textquotedblleft collapsing\textquotedblright%
\ $\mu_{n-1}^{\left(  c\right)  }\left[  \mathbf{g}\right]  =\mu_{n}\left[
\mathbf{g},c\right]  $, because $\mathbf{g}$ is a polyad of length $\left(
n-1\right)  $.
A universal algebra is a set which is closed under several polyadic operations
\cite{bergman2}.
If a concrete universal algebra has one fundamental
$n$-ary operation, called a \textit{polyadic multiplication} (or
$n$\textit{-ary multiplication}) $\mu_{n}$, we name it a ``polyadic system''\footnote{A set with one closed binary operation
without any other relations was called a groupoid by Hausmann and Ore
\cite{hau/ore} (see, also \cite{cli/pre1}). Nowadays the term
\textquotedblleft groupoid\textquotedblright\ is widely used in the category
theory and homotopy theory for a different construction, the so-called Brandt
groupoid \cite{bra1}. Bourbaki \cite{bourbaki1} introduced the term
\textquotedblleft magma\textquotedblright. To avoid misreading we will use the
neutral notation \textquotedblleft polyadic system\textquotedblright.}.  \begin{definition}
A \textit{polyadic system} $\mathfrak{G}=\left\langle \text{set}%
|\text{one fundamental operation}\right\rangle $ is a set $G$ which is closed under polyadic multiplication.
\end{definition}

More specifically, a $n$\textit{-ary system}
$\mathfrak{G}_{n}=\left\langle G\mid\mu_{n}\right\rangle $ is a set $G$ closed
under one $n$-ary operation $\mu_{n}$ (without any other additional structure).

For a given $n$-ary system $\left\langle G\mid\mu_{n}\right\rangle $ one can
construct another polyadic system $\left\langle G\mid\mu_{n^{\prime}}^{\prime
}\right\rangle $ over the same set $G$, but with another multiplication
$\mu_{n^{\prime}}^{\prime}$ of different arity $n^{\prime}$. In general, there
are three ways of changing the arity:

\begin{enumerate}
\item \textsl{Iterating}. Composition of the operation $\mu_{n}$ with itself
increases the arity from $n$ to $n^{\prime}=n_{iter}>n$. We denote the number
of iterating multiplications by $\ell_{\mu}$ and call the resulting
composition an \textit{iterated product}\footnote{Sometimes $\mathbf{\mu}%
_{n}^{\ell_{\mu}}$ is named a long product \cite{dor3}.} $\mathbf{\mu}%
_{n}^{\ell_{\mu}}$ (using the bold Greek letters) as (or $\mathbf{\mu}%
_{n}^{\bullet}$ if $\ell_{\mu}$ is obvious or not important)%
\begin{equation}
\mu_{n^{\prime}}^{\prime}=\mathbf{\mu}_{n}^{\ell_{\mu}}\overset{def}%
{=}\overset{\ell_{\mu}}{\overbrace{\mu_{n}\circ\left(  \mu_{n}\circ
\ldots\left(  \mu_{n}\times\operatorname*{id}\nolimits^{\times\left(
n-1\right)  }\right)  \ldots\times\operatorname*{id}\nolimits^{\times\left(
n-1\right)  }\right)  }},\label{mn}%
\end{equation}
where the final arity is%
\begin{equation}
n^{\prime}=n_{iter}=\ell_{\mu}\left(  n-1\right)  +1.\label{n}%
\end{equation}
There are many variants of placing $\mu_{n}$'s among $\operatorname*{id}${}'s
in the r.h.s. of (\ref{mn}), if no associativity is assumed.

An example of the iterated product can be given for a ternary operation
$\mu_{3}$ ($n=3$), where we can construct a 7-ary operation ($n^{\prime}=7$) by
$\ell_{\mu}=3$ compositions%
\begin{equation}
\mu_{7}^{\prime}\left[  g_{1},\ldots,g_{7}\right]  =\mathbf{\mu}_{3}%
^{3}\left[  g_{1},\ldots,g_{7}\right]  =\mu_{3}\left[  \mu_{3}\left[  \mu
_{3}\left[  g_{1},g_{2},g_{3}\right]  ,g_{4},g_{5}\right]  ,g_{6}%
,g_{7}\right]  , \label{m7}%
\end{equation}
and the corresponding commutative diagram is
\begin{equation}
\begin{diagram} G^{\times 7} & \rTo^{\mu_{3}\times\operatorname*{id}^{\times 4}} & G^{\times 5}&\rTo^{\mu_{3}\times\operatorname*{id}^{\times 2}}& G^{\times 3}\\ &\rdTo(5,3)_{\mu^{\prime}_{7}=\mathbf{\mu}_{3}^{3}} &&&\\ &&& & \dTo_{\mu_{3}} \\ && &&G \\ \end{diagram} \label{diam3}%
\end{equation}
In the general case, the horizontal part of the (iterating) diagram (\ref{diam3})
consists of $\ell_{\mu}$ terms. \medskip

\noindent

\item \textsl{Reducing (Collapsing)}. To decrease arity from $n$ to
$n^{\prime}=n_{red}<n$ one can use $n_{c}$ distinguished elements
(\textquotedblleft constants\textquotedblright) as additional $0$-ary
operations $\mu_{0}^{\left(  c_{i}\right)  }$, $i=1,\ldots n_{c}$, such
that\footnote{In \cite{dud/mic2} $\mu_{n}^{\left(  c_{1}\ldots c_{n_{c}%
}\right)  }$ is called a retract, which is already a busy and widely used term
in category theory for another construction.} the reduced product is defined
by%
\begin{equation}
\mu_{n^{\prime}}^{\prime}=\mu_{n^{\prime}}^{\left(  c_{1}\ldots c_{n_{c}%
}\right)  }\overset{def}{=}\mu_{n}\circ\left(  \overset{n_{c}}{\overbrace
{\mu_{0}^{\left(  c_{1}\right)  }\times\ldots\times\mu_{0}^{\left(  c_{n_{c}%
}\right)  }}}\times\operatorname*{id}\nolimits^{\times\left(  n-n_{c}\right)
}\right)  , \label{mc}%
\end{equation}
where%
\begin{equation}
n^{\prime}=n_{red}=n-n_{c}, \label{nr}%
\end{equation}
and the $0$-ary operations $\mu_{0}^{\left(  c_{i}\right)  }$ can be on any
places in (\ref{mc}). For instance, if we compose $\mu_{n}$ with the 0-ary
operation $\mu_{0}^{\left(  c\right)  }$, we obtain%
\begin{equation}
\mu_{n-1}^{\left(  c\right)  }\left[  \mathbf{g}\right]  =\mu_{n}\left[
\mathbf{g},c\right]  , \label{mcg}%
\end{equation}
and this reduced product is described by the commutative diagram
\begin{equation}
\begin{diagram} G^{\times \left( n-1 \right)}\times \left\{ \bullet \right\} & \rTo^{\operatorname*{id}^{\times \left( n-1 \right)}\times\mu_{0}^{\left( c\right) }} & G^{\times n} \\ \uTo^{\epsilon} & & \dTo_{\mu_{n}} \\ G^{\times \left( n-1 \right)} & \rTo^{ \mu_{n-1}^{\left( c\right)}} & G \\ \end{diagram} \label{dia0}%
\end{equation}
which can be treated as a definition of a new $\left(  n-1\right)  $-ary
operation $\mu_{n-1}^{\left(  c\right)  }=\mu_{n}\circ\mu_{0}^{\left(
c\right)  }$.

\medskip

\item \textsl{Mixing}. Changing (increasing or decreasing) arity by combining
the iterating and reducing (collapsing) methods.
\end{enumerate}

\begin{example}
If the initial multiplication is binary $\mu_{2}=\left(  \cdot\right)  $, and
there is one 0-ary operation $\mu_{0}^{\left(  c\right)  }$, we can construct
the following mixing operation%
\begin{equation}
\mu_{n}^{\left(  c\right)  }\left[  g_{1},\ldots,g_{n}\right]  =g_{1}\cdot
g_{2}\cdot\ldots\cdot g_{n}\cdot c, \label{mnc}%
\end{equation}
which in our notation can be called a $c$-iterated
multiplication\footnote{According to \cite{dud07} the operation (\ref{mnc})
can be called $c$-derived.}.
\end{example}

Let us recall some special elements of polyadic systems. A positive power of
an element (according to Post \cite{pos}) coincides with the number of
multiplications $\ell_{\mu}$ in the iteration (\ref{mn}).

\begin{definition}
A (positive) \textit{polyadic power} of an element is%
\begin{equation}
g^{\left\langle \ell_{\mu}\right\rangle }=\mathbf{\mu}_{n}^{\ell_{\mu}}\left[
g^{\ell_{\mu}\left(  n-1\right)  +1}\right]  . \label{pp}%
\end{equation}

\end{definition}

\begin{example}
\label{ex-qadd}Let us consider a polyadic version of the binary $q$-addition
which appears in study of nonextensive statistics (see, e.g.,
\cite{tsa94,niv/lem/wan})%
\begin{equation}
\mu_{n}\left[  \mathbf{g}\right]  =\overset{n}{\underset{i=1}{\sum}}%
g_{i}+\hbar\underset{i=1}{\overset{n}{%
{\displaystyle\prod}
}}g_{i}, \label{qadd}%
\end{equation}
where $g_{i}\in\mathbb{C}$ and $\hbar=1-q_{0}$, $q_{0}$ is a real constant (we
put here $q_{0}\neq1$ or $\hbar\neq0$). It is obvious that $g^{\left\langle
0\right\rangle }=g$, and%
\begin{equation}
g^{\left\langle 1\right\rangle }=\mu_{n}\left[  g^{n-1},g^{\left\langle
0\right\rangle }\right]  =ng+\hbar g^{n}.
\end{equation}
So we have the following recurrence formula%
\begin{equation}
g^{\left\langle k\right\rangle }=\mu_{n}\left[  g^{n-1},g^{\left\langle
k-1\right\rangle }\right]  =\left(  n-1\right)  g+\left(  1+\hbar
g^{n-1}\right)  g^{\left\langle k-1\right\rangle }.
\end{equation}
Solving this for an arbitrary polyadic power we get%
\begin{equation}
g^{\left\langle k\right\rangle }=g\left(  1+\dfrac{n-1}{\hbar}g^{1-n}\right)
\left(  1+\hbar g^{n-1}\right)  ^{k}-\dfrac{n-1}{\hbar}g^{2-n}.
\end{equation}

\end{example}

\begin{definition}
A \textit{polyadic (}$n$-\textit{ary) identity} (or neutral element) of a
polyadic system is a distinguished element $\varepsilon$ (and the
corresponding 0-ary operation $\mu_{0}^{\left(  \varepsilon\right)  }$) such
that for any element $g\in G$ we have \cite{rob58}%
\begin{equation}
\mu_{n}\left[  g,\varepsilon^{n-1}\right]  =g, \label{e}%
\end{equation}
where $g$ can be on any place in the l.h.s. of (\ref{e}).
\end{definition}

In polyadic systems, for an element $g$ there can exist many \textit{neutral
polyads} $\mathbf{n}\in G^{\times\left(  n-1\right)  }$ satisfying%
\begin{equation}
\mu_{n}\left[  g,\mathbf{n}\right]  =g,\label{mng}%
\end{equation}
where $g$ may be on any place. The neutral polyads are not determined
uniquely. It follows from (\ref{e}) and (\ref{mng}) that $\varepsilon^{n-1}$
is a neutral polyad.

\begin{definition}
\label{def-mi}An element of a polyadic system $g$ is called $\ell_{\mu}%
$-\textit{idempotent}, if there exist such $\ell_{\mu}$ that%
\begin{equation}
g^{\left\langle \ell_{\mu}\right\rangle }=g. \label{mi}%
\end{equation}

\end{definition}

It is obvious that an identity is $\ell_{\mu}$-idempotent with arbitrary
$\ell_{\mu}$. We define \textit{(total) associativity} as invariance of the
composition of two $n$-ary multiplications%
\begin{equation}
\mathbf{\mu}_{n}^{2}\left[  \mathbf{g},\mathbf{h},\mathbf{u}\right]
=invariant\label{ghu}%
\end{equation}
under placement of the internal multiplication in the r.h.s. with a fixed order
of elements in the whole polyad of $\left(  2n-1\right)  $ elements
$\mathbf{t}^{\left(  2n-1\right)  }=\left(  \mathbf{g},\mathbf{h}%
,\mathbf{u}\right)  $. Informally, \textquotedblleft internal
brackets/multiplication can be moved on any place\textquotedblright, which
gives%
\begin{equation}
\mu_{n}\circ\left(  \overset{i=1}{\mu_{n}}\times\operatorname*{id}%
\nolimits^{\times\left(  n-1\right)  }\right)  =\mu_{n}\circ\left(
\operatorname*{id}\times\overset{i=2}{\mu_{n}}\times\operatorname*{id}%
\nolimits^{\times\left(  n-2\right)  }\right)  =\ldots=\mu_{n}\circ\left(
\operatorname*{id}\nolimits^{\times\left(  n-1\right)  }\times\overset
{i=n}{\mu_{n}}\right)  ,
\end{equation}
where the internal $\mu_{n}$ can be on any place $i=1,\ldots,n$. There are
many other particular kinds of associativity which were introduced in
\cite{pos,thu49} and studied in \cite{sus29,sok1} (see, also \cite{cou10}).
Here we will confine ourselves to the most general, total associativity
(\ref{ghu}).

\begin{definition}
A \textit{polyadic semigroup} ($n$-\textit{ary semigroup}) is a $n$-ary system
whose operation is associative, or $\mathfrak{G}_{n}^{semigrp}=\left\langle
G\mid\mu_{n}\mid\text{associativity (\ref{ghu})}\right\rangle $.
\end{definition}

In general, it is very important to find the \textit{associativity preserving}
\textit{conditions}, when an associative initial operation $\mu_{n}$ leads to
an associative final operation $\mu_{n^{\prime}}^{\prime}$ while changing the
arity (by iterating (\ref{mn}) or reducing (\ref{mc})).

\begin{example}
An associativity preserving reduction can be given by the construction of a
binary associative operation using a $\left(  n-2\right)  $-tuple $\mathbf{c}$
as%
\begin{equation}
\mu_{2}^{\left(  \mathbf{c}\right)  }\left[  g,h\right]  =\mu_{n}\left[
g,\mathbf{c},h\right]  . \label{mgh}%
\end{equation}

\end{example}

The associativity preserving mixing constructions with different arities and
places were considered in \cite{dud/mic2,mic2,sok1}.

In polyadic systems, there are several analogs of binary commutativity. The
most straightforward one comes from commutation of the multiplication with permutations.

\begin{definition}
A polyadic system is $\sigma$-\textit{commutative}, if $\mu_{n}=\mu_{n}%
\circ\sigma$, where $\sigma$ is a fixed element of $S_{n}$, the permutation
group on $n$ elements. If this holds for all $\sigma\in S_{n}$, then a
polyadic system is \textit{commutative}.
\end{definition}

A special type of the $\sigma$-commutativity%
\begin{equation}
\mu_{n}\left[  g,\mathbf{t},h\right]  =\mu_{n}\left[  h,\mathbf{t},g\right]
\label{mth}%
\end{equation}
is called \textit{semicommutativity}. So for a $n$-ary semicommutative system
we have%
\begin{equation}
\mu_{n}\left[  g,h^{n-1}\right]  =\mu_{n}\left[  h^{n-1},g\right]  .
\end{equation}

If a $n$-ary semigroup $\mathfrak{G}_{n}^{semigrp}$ is iterated from a
commutative binary semigroup with identity, then $\mathfrak{G}_{n}^{semigrp}$
is semicommutative. Another possibility is to generalize the binary mediality
in semigroups%
\begin{equation}
\left(  g_{11}\cdot g_{12}\right)  \cdot\left(  g_{21}\cdot g_{22}\right)
=\left(  g_{11}\cdot g_{21}\right)  \cdot\left(  g_{12}\cdot g_{22}\right)  ,
\end{equation}
which follows from the binary commutativity. For $n$-ary systems, it is seen
that the mediality should contain $\left(  n+1\right)  $ multiplications, that
it is a relation between $n\times n$ elements, and therefore that it can be presented
in a matrix from.

\begin{definition}
A polyadic system is \textit{medial} (or \textit{entropic}), if
\cite{eva7,belousov}%
\begin{equation}
\mu_{n}\left[
\begin{array}
[c]{c}%
\mu_{n}\left[  g_{11},\ldots,g_{1n}\right] \\
\vdots\\
\mu_{n}\left[  g_{n1},\ldots,g_{nn}\right]
\end{array}
\right]  =\mu_{n}\left[
\begin{array}
[c]{c}%
\mu_{n}\left[  g_{11},\ldots,g_{n1}\right] \\
\vdots\\
\mu_{n}\left[  g_{1n},\ldots,g_{nn}\right]
\end{array}
\right]  . \label{mmn}%
\end{equation}

\end{definition}

In the case of polyadic semigroups we use the notation (\ref{mn}) and can
present the mediality as follows%
\begin{equation}
\mathbf{\mu}_{n}^{n}\left[  \mathbf{G}\right]  =\mathbf{\mu}_{n}^{n}\left[
\mathbf{G}^{T}\right]  ,
\end{equation}
where $\mathbf{G}=\left\Vert g_{ij}\right\Vert $ is the $n\times n$ matrix of
elements and $\mathbf{G}^{T}$ is its transpose.

The semicommutative polyadic semigroups are medial, as in the binary case,
but, in general (except $n=3$) not vice versa \cite{gla/gle}.

\begin{definition}
\label{def-cancel}A polyadic system is \textit{cancellative}, if%
\begin{equation}
\mu_{n}\left[  g,\mathbf{t}\right]  =\mu_{n}\left[  h,\mathbf{t}\right]
\Longrightarrow g=h, \label{mg}%
\end{equation}
where $g,h$ can be on any place. This means that the mapping $\mu_{n}$ is
one-to-one in each variable.
\end{definition}

If $g,h$ are on the same $i$-th place on both sides of (\ref{mg}), the
polyadic system is called $i$-\textit{cancellative}. The \textit{left} and
\textit{right} cancellativity are $1$-cancellativity and $n$-cancellativity
respectively. A right and left cancellative $n$-ary semigroup is cancellative
(with respect to the same subset).

\begin{definition}
A polyadic system is called (uniquely) $i$-\textit{solvable}, if for all
polyads $\mathbf{t}$, $\mathbf{u}$ and element $h$, one can (uniquely) resolve
the equation (with respect to $h$) for the fundamental operation%
\begin{equation}
\mu_{n}\left[  \mathbf{u},h,\mathbf{t}\right]  =g \label{mug}%
\end{equation}
where $h$ can be on any $i$-th place.
\end{definition}

\begin{definition}
A polyadic system which is uniquely $i$-solvable for all places $i=1,\ldots,n$
in (\ref{mug}) is called a $n$-\textit{ary }(or \textit{polyadic})\textit{
quasigroup}.
\end{definition}

It follows, that, if (\ref{mug}) uniquely $i$-solvable for all places, then%
\begin{equation}
\mathbf{\mu}_{n}^{\ell_{\mu}}\left[  \mathbf{u},h,\mathbf{t}\right]  =g
\label{mugl}%
\end{equation}
can be (uniquely) resolved with respect to $h$ being on any place.

\begin{definition}
\label{def-grp}An associative polyadic quasigroup is called a $n$-\textit{ary}
(or \textit{polyadic})\textit{ group}.
\end{definition}

In a polyadic group the only solution of (\ref{mug}) is called a
\textit{querelement}\footnote{We use the original notation after \cite{dor3}
and do not use \textquotedblleft skew element\textquotedblright, because
it\ can be confused with the wide usage of \textquotedblleft
skew\textquotedblright\ in other, different senses.} of $g$ and is denoted by
$\bar{g}$ \cite{dor3}, such that%
\begin{equation}
\mu_{n}\left[  \mathbf{h},\bar{g}\right]  =g, \label{mgg}%
\end{equation}
where $\bar{g}$ can be on any place. Obviously, any idempotent $g$ coincides
with its querelement $\bar{g}=g$.

\begin{example}
For the $q$-addition (\ref{qadd}) from Example \ref{ex-qadd}, using
(\ref{mgg}) with $\mathbf{h}=g^{n-1}$ we obtain%
\begin{equation}
\bar{g}=-\dfrac{\left(  n-2\right)  g}{1+\hbar g^{n-1}}.
\end{equation}

\end{example}

It follows from (\ref{mgg}) and (\ref{mng}), that the polyad%
\begin{equation}
\mathbf{n}_{\left(  \bar{g}\right)  }=\left(  g^{n-2},\bar{g}\right)
\label{ng}%
\end{equation}
is neutral for any element $g$, where $\bar{g}$ can be on any place. If this
$i$-th place is important, then we write $\mathbf{n}_{\left(  g\right)  ,i}$.
More generally, because any neutral polyad plays a role of identity (see
(\ref{mng})), for any element $g$ we define its \textit{polyadic inverse} (the
sequence of length $\left(  n-2\right)  $ denoted by the same letter
$\mathbf{g}^{-1}$ in bold) as (see \cite{pos} and by modified analogy with
\cite{pet05,galmak1})%
\begin{equation}
\mathbf{n}_{\left(  g\right)  }=\left(  \mathbf{g}^{-1},g\right)  =\left(
g,\mathbf{g}^{-1}\right)  ,\label{ngg}%
\end{equation}
which can be written in terms of the multiplication as%
\begin{equation}
\mu_{n}\left[  g,\mathbf{g}^{-1},h\right]  =\mu_{n}\left[  h,\mathbf{g}%
^{-1},g\right]  =h\label{mng1}%
\end{equation}
for all $h$ in $G$. It is obvious that the polyads
\begin{equation}
\mathbf{n}_{\left(  g^{k}\right)  }=\left(  \left(  \mathbf{g}^{-1}\right)
^{k},g^{k}\right)  =\left(  g^{k},\left(  \mathbf{g}^{-1}\right)  ^{k}\right)
\label{ngk}%
\end{equation}
are neutral as well for any $k\geq1$. It follows from (\ref{ng}) that the
polyadic inverse of $g$ is $\left(  g^{n-3},\bar{g}\right)  $, and one of
$\bar{g}$ is $\left(  g^{n-2}\right)  $, and in this case $g$ is called
\textit{querable}. In a polyadic group all elements are querable
\cite{dud/gla/gle,cel1}.

The number of relations in (\ref{mgg}) can be reduced from $n$ (the number of
possible places) to only $2$ (when $g$ is on the first and last places
\cite{dor3,tim72}), such that in a polyadic group the \textit{D\"{o}rnte
relations}%
\begin{equation}
\mu_{n}\left[  g,\mathbf{n}_{\left(  g\right)  ,i}\right]  =\mu_{n}\left[
\mathbf{n}_{\left(  g\right)  ,j},g\right]  =g \label{mgnn}%
\end{equation}
hold valid for any allowable $i,j$, and (\ref{mgnn}) are analogs of $g\cdot
h\cdot h^{-1}=h\cdot h^{-1}\cdot g=g$ in binary groups. The relation
(\ref{mgg}) can be treated as a definition of the (unary)
\textit{queroperation} $\bar{\mu}_{1}:G\rightarrow G$ by%
\begin{equation}
\bar{\mu}_{1}\left[  g\right]  =\bar{g}, \label{m1g}%
\end{equation}
such that the diagram%
\begin{equation}
\begin{diagram} G^{\times n} &\rTo^{\mu_{n}} & G\\ \uTo^{\operatorname*{id}^{\times \left( n-1 \right)}\times \bar{\mu}_{1}}&\;\;\ruTo(3,2)_{\mathsf{Pr}_n}&\\ G^{\times n}& & \end{diagram} \label{diagm4}%
\end{equation}

\medskip

\noindent commutes. Then, using the queroperation (\ref{m1g}) one can give a
diagrammatic definition of a polyadic group (cf. \cite{gle/gla}).

\begin{definition}
A \textit{polyadic group} is a universal algebra
\begin{equation}
\mathfrak{G}_{n}^{grp}=\left\langle G\mid\mu_{n},\bar{\mu}_{1}\mid
\text{associativity, D\"{o}rnte relations}\right\rangle ,
\end{equation}
where $\mu_{n}$ is $n$-ary associative operation and $\bar{\mu}_{1}$ is the
queroperation (\ref{m1g}), such that the following diagram%
\begin{equation}
\begin{diagram} G^{\times \left( n \right)} & \rTo^{\operatorname*{id}^{\times \left( n-1 \right)} \times\bar{\mu}_{1}}& & G^{\times n}&& \lTo^{\;\;\bar{\mu}_{1}\times\operatorname*{id}^{\times \left( n-1 \right)} } & G^{\times n} \\ \uTo^{\operatorname*{id}\times\mathsf{Diag}_{ \left( n-1 \right)}} && & \dTo_{\mu_{n}}&& & \uTo_{\mathsf{Diag}_{ \left( n-1 \right)}\times\operatorname*{id}} \\ G\times G & \rTo^{ \mathsf{Pr}_1} && G& &\lTo^{ \mathsf{Pr}_2} & G\times G\\ \end{diagram} \label{diam5}%
\end{equation}
\medskip

\noindent commutes, where $\bar{\mu}_{1}$ can be only on the first and second
places from the right (resp. left) on the left (resp. right) part of the diagram.
\end{definition}

A straightforward generalization of the queroperation concept and
corresponding definitions can be made by substituting in the above formulas
(\ref{mgg})--(\ref{m1g}) the $n$-ary multiplication $\mu_{n}$ by the iterating
multiplication $\mathbf{\mu}_{n}^{\ell_{\mu}}$ (\ref{mn}) (cf. \cite{dud2} for
$\ell_{\mu}=2$ and \cite{galmak2}).

Let us define the \textit{querpower} $k$ of $g$ recursively by
\cite{wan/wan79,mic83}%
\begin{equation}
\bar{g}^{\left\langle \left\langle k\right\rangle \right\rangle }%
=\overline{\left(  \bar{g}^{\left\langle \left\langle k-1\right\rangle
\right\rangle }\right)  },\label{gk}%
\end{equation}
where $\bar{g}^{\left\langle \left\langle 0\right\rangle \right\rangle }=g$,
$\bar{g}^{\left\langle \left\langle 1\right\rangle \right\rangle }=\bar{g}$,
$\bar{g}^{\left\langle \left\langle 2\right\rangle \right\rangle }%
=\overline{\bar{g}}$,... or as the $k$ composition $\bar{\mu}_{1}^{\circ
k}=\overset{k}{\overbrace{\bar{\mu}_{1}\circ\bar{\mu}_{1}\circ\ldots\circ
\bar{\mu}_{1}}}$ of the unary queroperation (\ref{m1g}). We can define the
\textit{negative polyadic power} of an element $g$ by the recursive
relationship%
\begin{equation}
\mu_{n}\left[  g^{\left\langle \ell_{\mu}-1\right\rangle },g^{n-2}%
,g^{\left\langle -\ell_{\mu}\right\rangle }\right]  =g,
\end{equation}
or (after the use of the positive polyadic power (\ref{pp})) as a solution of
the equation%
\begin{equation}
\mathbf{\mu}_{n}^{\ell_{\mu}}\left[  g^{\ell_{\mu}\left(  n-1\right)
},g^{\left\langle -\ell_{\mu}\right\rangle }\right]  =g.\label{ml}%
\end{equation}

The querpower (\ref{gk}) and the polyadic power (\ref{ml}) are connected
\cite{dud1}. We reformulate this connection using the so called Heine numbers
\cite{heine} or $q$-deformed numbers \cite{kac/che}%
\begin{equation}
\left[  \left[  k\right]  \right]  _{q}=\dfrac{q^{k}-1}{q-1}, \label{kq}%
\end{equation}
which have the \textquotedblleft nondeformed\textquotedblright\ limit
$q\rightarrow1$ as $\left[  \left[  k\right]  \right]  _{q}\rightarrow k$ and
$\left[  \left[  0\right]  \right]  _{q}=0$. If $\left[  \left[  k\right]
\right]  _{q}=0$, then $q$ is a $k$-th root of unity. From (\ref{gk}) and
(\ref{ml}) we obtain%
\begin{equation}
\bar{g}^{\left\langle \left\langle k\right\rangle \right\rangle }%
=g^{\left\langle -\left[  \left[  k\right]  \right]  _{2-n}\right\rangle },
\end{equation}
which can be treated as the following \textquotedblleft
deformation\textquotedblright\ statement:

\begin{assertion}
The querpower coincides with the negative polyadic deformed\textsl{ }power
with the \textquotedblleft deformation\textquotedblright\ parameter $q$ which
is equal to the \textquotedblleft deviation\textquotedblright\ $\left(
2-n\right)  $ from the binary group.
\end{assertion}

\begin{example}
Let us consider a binary group $\mathfrak{G}_{2}=\left\langle G\mid\mu
_{2}\right\rangle $, we denote $\mu_{2}=\left(  \cdot\right)  $, and construct
(using (\ref{mn}) and (\ref{mc})) the reduced $4$-ary product by $\mu
_{4}^{\prime}\left[  \mathbf{g}\right]  =g_{1}\cdot g_{2}\cdot g_{3}\cdot
g_{4}\cdot c$, where $g_{i}\in G$ and $c$ is in the center of the group
$\mathfrak{G}_{2}$. In the 4-ary group $\mathfrak{G}_{4}^{\prime}=\left\langle
G,\mu_{4}^{\prime}\right\rangle $ we derive the following positive and
negative polyadic powers (obviously $g^{\left\langle 0\right\rangle }=\bar
{g}^{\left\langle \left\langle 0\right\rangle \right\rangle }=g$)%
\begin{align}
\ \ g^{\left\langle 1\right\rangle }  &  =g^{4}\cdot c,\ \ g^{\left\langle
2\right\rangle }=g^{7}\cdot c^{2},\ldots,g^{\left\langle k\right\rangle
}=g^{3k+1}\cdot c^{k},\\
\ \ g^{\left\langle -1\right\rangle }  &  =g^{-2}\cdot c^{-1}%
,\ \ g^{\left\langle -2\right\rangle }=g^{-5}\cdot c^{-2},\ldots
,g^{\left\langle -k\right\rangle }=g^{-3k+1}\cdot c^{-k},
\end{align}
and the querpowers%
\begin{equation}
\bar{g}^{\left\langle \left\langle 1\right\rangle \right\rangle }=g^{-2}\cdot
c^{-1},\ \ \bar{g}^{\left\langle \left\langle 2\right\rangle \right\rangle
}=g^{-4}\cdot c,\ldots,\ \ \bar{g}^{\left\langle \left\langle k\right\rangle
\right\rangle }=g^{\left(  -2\right)  ^{k}}\cdot c^{\left[  \left[  k\right]
\right]  _{-2}}.
\end{equation}

\end{example}

Let $\mathfrak{G}$$_{n}=\left\langle G\mid\mu_{n}\right\rangle $ and
$\mathfrak{G}$$_{n^{\prime}}^{\prime}=\left\langle G^{\prime}\mid
\mu_{n^{\prime}}^{\prime}\right\rangle $ be two polyadic systems of any kind.
If their multiplications are of the same arity $n=n^{\prime}$, then one can
define the following \textit{one-place} mappings from $\mathfrak{G}$$_{n}$ to
$\mathfrak{G}$$_{n}^{\prime}$ (for \textit{many-place} mappings, which
\textit{change} arity $n\neq n^{\prime}$ and corresonding
\textit{heteromorphisms}, see \cite{dup2012}).

Suppose we have $n+1$ mappings $\Phi_{i}:G\rightarrow G^{\prime}$,
$i=1,\ldots,n+1$. An ordered system of mappings $\left\{  \Phi_{i}\right\}  $
is called a \textit{homotopy} from $\mathfrak{G}$$_{n}$ to $\mathfrak{G}$%
$_{n}^{\prime}$, if (see, e.g., \cite{belousov})%
\begin{equation}
\Phi_{n+1}\left(  \mu_{n}\left[  g_{1},\ldots,g_{n}\right]  \right)  =\mu
_{n}^{\prime}\left[  \Phi_{1}\left(  g_{1}\right)  ,\ldots,\Phi_{n}\left(
g_{n}\right)  \right]  ,\ \ \ \ \ g_{i}\in G.\label{fm}%
\end{equation}

A \textit{homomorphism} from $\mathfrak{G}$$_{n}$ to $\mathfrak{G}$%
$_{n}^{\prime}$ is given, if there exists a (one-place) mapping $\Phi
:G\rightarrow G^{\prime}$ satisfying
\begin{equation}
\Phi\left(  \mu_{n}\left[  g_{1},\ldots,g_{n}\right]  \right)  =\mu
_{n}^{\prime}\left[  \Phi\left(  g_{1}\right)  ,\ldots,\Phi\left(
g_{n}\right)  \right]  ,\ \ \ \ \ g_{i}\in G, \label{fm1}%
\end{equation}
which means that the corresponding (\textit{equiary}\footnote{The map is
equiary, if it does not change the arity of operations i.e. $n=n^{\prime}$, for
nonequiary maps see \cite{dup2012} and refs. therein.}) diagram is commutative%
\begin{equation}
\begin{diagram} G & \rTo^{\varphi} & G^{\prime} \\ \uTo^{\mu_{n}} & & \uTo_{\mu_{n}^{\prime}} \\ G^{\times n} & \rTo^{ \left( \varphi \right) ^{\times n}} & \left( G^{\prime}\right)^{\times n} \\ \end{diagram} \label{dia2}%
\end{equation}
\medskip It is obvious that, if a polyadic system contains distinguished
elements (identities, querelements, etc.), they are also mapped by $\varphi$
correspondingly (for details and a review, see, e.g., \cite{galmak2,rusakov1}).
The most important application of one-place mappings is in establishing a
general structure for $n$-ary multiplication.

\section{The Hossz\'{u}-Gluskin theorem}

Let us consider possible concrete forms of polyadic multiplication in terms of
lesser arity operations. Obviously, the simplest way of constructing a $n$-ary
product $\mu_{n}^{\prime}$ from the binary one $\mu_{2}=\left(  \ast\right)  $
is $\ell_{\mu}=n$ iteration (\ref{mn}) \cite{sus35,mil35}%
\begin{equation}
\mu_{n}^{\prime}\left[  \mathbf{g}\right]  =g_{1}\ast g_{2}\ast\ldots\ast
g_{n},\ \ \ \ g_{i}\in G.\label{mgn}%
\end{equation}
In \cite{dor3} it was noted that not all $n$-ary groups have a product of this
special form. The binary group $\mathfrak{G}_{2}^{\ast}=\left\langle G\mid
\mu_{2}=\ast,e\right\rangle $ was called a \textit{covering group} of the
$n$-ary group $\mathfrak{G}_{n}^{\prime}=\left\langle G\mid\mu_{n}^{\prime
}\right\rangle $ in \cite{pos} (see, also, \cite{tve53}), where a theorem
establishing a more general (than (\ref{mgn})) structure of $\mu_{n}^{\prime
}\left[  \mathbf{g}\right]  $ in terms of subgroup structure of the covering
group was given. A manifest form of the $n$-ary group product $\mu_{n}%
^{\prime}\left[  \mathbf{g}\right]  $ in terms of the binary one and a special
mapping was found in \cite{hos,glu1} and is called the Hossz\'{u}-Gluskin
theorem, despite the same formulas having appeared much earlier in \cite{tur38,pos}
(for the relationship between the formulations, see \cite{gal/vor}). A simple
construction of $\mu_{n}^{\prime}\left[  \mathbf{g}\right]  $ which is present
in the Hossz\'{u}-Gluskin theorem was given in \cite{sok}. Here we follow this
scheme in the opposite direction, by just deriving the final formula step by step
(without writing it immediately) with clear examples. Then we introduce a
\textquotedblleft deformation\textquotedblright\ to it in such a way that a
generalized \textquotedblleft$q$-deformed\textquotedblright%
\ Hossz\'{u}-Gluskin theorem can be formulated.

First, let us rewrite (\ref{mgn}) in its equivalent form%

\begin{equation}
\mu_{n}^{\prime}\left[  \mathbf{g}\right]  =g_{1}\ast g_{2}\ast\ldots\ast
g_{n}\ast e,\ \ \ \ g_{i},e\in G, \label{me}%
\end{equation}
where $e$ is a distinguished element of the binary group $\left\langle
G\mid\ast,e\right\rangle $, that is the identity. Now we apply to (\ref{me})
an \textquotedblleft extended\textquotedblright\ version of the homotopy
relation (\ref{fm}) with $\Phi_{i}=\psi_{i}$, $i=1,\ldots n$, and the l.h.s.
mapping $\Phi_{n+1}=\operatorname*{id}$, but add an action $\psi_{n+1}$ on the
identity $e$ of the binary group $\left\langle G\mid\ast,e\right\rangle $.
Then we get (see (\ref{mcg}) and (\ref{mnc}))%
\begin{equation}
\mu_{n}\left[  \mathbf{g}\right]  =\mu_{n}^{\left(  e\right)  }\left[
\mathbf{g}\right]  =\psi_{1}\left(  g_{1}\right)  \ast\psi_{2}\left(
g_{2}\right)  \ast\ldots\ast\psi_{n}\left(  g_{n}\right)  \ast\psi
_{n+1}\left(  e\right)  =\left(  \ast\underset{i=1}{\overset{n}{%
{\displaystyle\prod}
}}\psi_{i}\left(  g_{i}\right)  \right)  \ast\psi_{n+1}\left(  e\right)  .
\label{mge}%
\end{equation}

In this way we have obtained the most general form of polyadic multiplication
in terms of $\left(  n+1\right)  $ \textquotedblleft
extended\textquotedblright\ homotopy maps $\psi_{i}$, $i=1,\ldots n+1$, such
that the diagram
\begin{equation}
\begin{diagram} G^{\times \left( n \right)}\times \left\{ \bullet \right\} & \rTo^{\operatorname*{id}^{\times n }\times\mu_{0}^{\left( e\right) }} & G^{\times \left( n+1 \right)} & \rTo^{\psi_{1} \times \ldots \times \psi_{n+1}} & G^{\times \left( n+1 \right)} \\ \uTo^{\epsilon} & & & & \dTo_{\mu_{2}^{\times n}} \\ G^{\times \left( n \right)} & & \rTo^{ \mu_{n}^{\left( e\right)}} & & G \\ \end{diagram}
\end{equation}
\noindent commutes. A natural question arises, whether all associative
polyadic systems have this form of multiplication or do we have others? In
general, we can correspondingly classify polyadic systems as:%
\begin{align}
&  \text{1) \textit{Homotopic }polyadic systems which can be presented in the
form (\ref{mge}).}\label{1}\\
&  \text{2) \textit{Nonhomotopic} polyadic systems with multiplication of
other than (\ref{mge}) shapes. }\label{2}%
\end{align}
If the second class is nonempty, it would be interesting to find examples of
nonhomotopic polyadic systems. The Hossz\'{u}-Gluskin theorem
considers the homotopic polyadic systems and gives one of the possible choices for
the \textquotedblleft extended\textquotedblright\ homotopy maps $\psi_{i}$ in
(\ref{mge}). We will show that this choice can be extended (\textquotedblleft
deformed\textquotedblright) to the infinite \textquotedblleft
q-series\textquotedblright.

The main idea in constructing the \textquotedblleft
automatically\textquotedblright\ associative $n$-ary operation $\mu_{n}$ in
(\ref{mge}) is to express the binary multiplication $\left(  \ast\right)  $
and the \textquotedblleft extended\textquotedblright\ homotopy maps $\psi_{i}$
in terms of $\mu_{n}$ itself \cite{sok}. A simplest binary multiplication
which can be built from $\mu_{n}$ is (see (\ref{mgh}))%
\begin{equation}
g\ast_{\mathbf{t}}h=\mu_{n}\left[  g,\mathbf{t},h\right]  ,
\end{equation}
where $\mathbf{t}$ is any fixed polyad of length $\left(  n-2\right)  $. If we
apply here the equations for the identity $e$ in a binary group%
\begin{equation}
g\ast_{\mathbf{t}}e=g,\ \ \ \ e\ast_{\mathbf{t}}h=h,\label{ge}%
\end{equation}
then we obtain%
\begin{equation}
\mu_{n}\left[  g,\mathbf{t},e\right]  =g,\ \ \ \mu_{n}\left[  e,\mathbf{t}%
,h\right]  =h.\label{mt}%
\end{equation}

We observe from (\ref{mt}) that $\left(  \mathbf{t},e\right)  $ and $\left(
e,\mathbf{t}\right)  $ are neutral sequences of length $\left(  n-1\right)  $,
and therefore using (\ref{ngg}) we can take $\mathbf{t}$ as a polyadic inverse
of $e$ (the identity of the binary group) considered as an element (but not an
identity) of the polyadic system $\left\langle G\mid\mu_{n}\right\rangle $,
that is $\mathbf{t}=\mathbf{e}^{-1}$. Then, the binary multiplication
constructed from $\mu_{n}$ and which has the standard identity properties
(\ref{ge}) can be chosen as%
\begin{equation}
g\ast h=g\ast_{e}h=\mu_{n}\left[  g,\mathbf{e}^{-1},h\right]  . \label{ghe}%
\end{equation}
Using this construction any element of the polyadic system $\left\langle G\mid\mu_{n}\right \rangle $ can be distinguished and may serve as the identity of the binary group, and is then denoted by $e$ (for clarity and convenience).  

We recognize in (\ref{ghe}) a version of the Maltsev term (see, e.g., \cite{bergman2}),
which can be called a \textit{polyadic Maltsev term} and is defined as%
\begin{equation}
p\left(  g,e,h\right)  \overset{def}{=}\mu_{n}\left[  g,\mathbf{e}%
^{-1},h\right]
\end{equation}
having the standard term properties \cite{bergman2}
\begin{equation}
p\left(  g,e,e\right)  =g,\ \ \ p\left(  e,e,h\right)  =h,
\end{equation}
which now follow from (\ref{mt}), i.e. the polyads $\left(  e,\mathbf{e}%
^{-1}\right)  $ and $\left(  \mathbf{e}^{-1},e\right)  $ are neutral, as they
should be (\ref{ngg}). Denote by $g^{-1}$ the inverse element of $g$ in the
binary group ($g\ast g^{-1}=g^{-1}\ast g=e$) and $\mathbf{g}^{-1}$ its
polyadic inverse in a $n$-ary group (\ref{ngg}), then it follows from
(\ref{ghe}) that $\mu_{n}\left[  g,\mathbf{e}^{-1},g^{-1}\right]  =e$. Thus,
we get%
\begin{equation}
g^{-1}=\mu_{n}\left[  e,\mathbf{g}^{-1},e\right]  , \label{mne}%
\end{equation}
which can be considered as a connection between the inverse $g^{-1}$ in the
binary group and the polyadic inverse in the polyadic system related to the
same element $g$. For $n$-ary group we can write $\mathbf{g}^{-1}=\left(
g^{n-3},\bar{g}\right)  $ and the binary group inverse $g^{-1}$ becomes%
\begin{equation}
g^{-1}=\mu_{n}\left[  e,g^{n-3},\bar{g},e\right]  .
\end{equation}

If $\left\langle G\mid\mu_{n}\right\rangle $ is a $n$-ary group, then the
element $e$ is querable (\ref{mng1}), for the polyadic inverse $\mathbf{e}%
^{-1}$ one can choose $\left(  e^{n-3},\bar{e}\right)  $ with $\bar{e}$ being
on any place, and the polyadic Maltsev term becomes \cite{shc03} $p\left(
g,e,h\right)  =\mu_{n}\left[  g,e^{n-3},\bar{e},h\right]  $ (together with the
multiplication (\ref{ghe})). For instance, if $n=3$, we have%
\begin{equation}
g\ast h=\mu_{3}\left[  g,\bar{e},h\right]  ,\ \ g^{-1}=\mu_{3}\left[
e,\bar{g},e\right]  , \label{g3}%
\end{equation}
and the neutral polyads are $\left(  e,\bar{e}\right)  $ and $\left(  \bar
{e},e\right)  $.

Now let us turn to build the main construction, that of the \textquotedblleft
extended\textquotedblright\ homotopy maps $\psi_{i}$ (\ref{mge}) in terms of
$\mu_{n}$, which will lead to the Hossz\'{u}-Gluskin theorem. We start with a
simple example of a ternary system (\ref{g3}), derive the Hossz\'{u}-Gluskin
\textquotedblleft chain formula\textquotedblright, and then it will be clear
how to proceed for generic $n$. Instead of (\ref{mge}) we write%
\begin{equation}
\mu_{3}\left[  g,h,u\right]  =\psi_{1}\left(  g\right)  \ast\psi_{2}\left(
h\right)  \ast\psi_{3}\left(  u\right)  \ast\psi_{4}\left(  e\right)
\label{m3}%
\end{equation}
and try to construct $\psi_{i}$ in terms of the ternary product $\mu_{3}$ and the
binary identity $e$. We already know the structure of the binary
multiplication (\ref{g3}): it contains $\bar{e}$, and therefore we can insert
between $g$, $h$ and $u$ in the l.h.s. of (\ref{m3}) a neutral ternary polyad
$\left(  \bar{e},e\right)  $ or its powers $\left(  \bar{e}^{k},e^{k}\right)
$. Thus, taking for all insertions the \textit{minimal number} of neutral
polyads, we get%
\begin{align}
\mu_{3}\left[  g,h,u\right]   &  =\mathbf{\mu}_{3}^{2}\left[  g,\overset{%
\begin{array}
[c]{c}%
\ast\\
\downarrow
\end{array}
}{\bar{e}},e,h,u\right]  =\mathbf{\mu}_{3}^{4}\left[  g,\overset{%
\begin{array}
[c]{c}%
\ast\\
\downarrow
\end{array}
}{\bar{e}},e,h,\bar{e},\overset{%
\begin{array}
[c]{c}%
\ast\\
\downarrow
\end{array}
}{\bar{e}},e,e,u\right] \nonumber\\
&  =\mathbf{\mu}_{3}^{7}\left[  g,\overset{%
\begin{array}
[c]{c}%
\ast\\
\downarrow
\end{array}
}{\bar{e}},e,h,\bar{e},\overset{%
\begin{array}
[c]{c}%
\ast\\
\downarrow
\end{array}
}{\bar{e}},e,e,u,\bar{e},\bar{e},\overset{%
\begin{array}
[c]{c}%
\ast\\
\downarrow
\end{array}
}{\bar{e}},e,e,e\right]  . \label{m3g}%
\end{align}

We show by arrows the binary products in special places: there should be
$1,3,5,\ldots\left(  2k-1\right)  $ elements in between them to form inner
ternary products. Then we rewrite (\ref{m3g}) as%
\begin{equation}
\mu_{3}\left[  g,h,u\right]  =\mathbf{\mu}_{3}^{3}\left[  g,\overset{%
\begin{array}
[c]{c}%
\ast\\
\downarrow
\end{array}
}{\bar{e}},\mu_{3}\left[  e,h,\bar{e}\right]  ,\overset{%
\begin{array}
[c]{c}%
\ast\\
\downarrow
\end{array}
}{\bar{e}},\mathbf{\mu}_{3}^{2}\left[  e,e,u,\bar{e},\bar{e}\right]
,\overset{%
\begin{array}
[c]{c}%
\ast\\
\downarrow
\end{array}
}{\bar{e}},\mu_{3}\left[  e,e,e\right]  \right]  . \label{m3h}%
\end{equation}
Comparing this with (\ref{m3}), we can exactly identify the \textquotedblleft
extended\textquotedblright\ homotopy maps $\psi_{i}$ as%
\begin{align}
\psi_{1}\left(  g\right)   &  =g,\label{p1}\\
\psi_{2}\left(  g\right)   &  =\varphi\left(  g\right)  ,\\
\psi_{3}\left(  g\right)   &  =\varphi\left(  \varphi\left(  g\right)
\right)  =\varphi^{2}\left(  g\right)  ,\\
\psi_{4}\left(  e\right)   &  =\mu_{3}\left[  e,e,e\right]  , \label{p4}%
\end{align}
where%
\begin{equation}
\varphi\left(  g\right)  =\mu_{3}\left[  e,g,\bar{e}\right]  , \label{fg}%
\end{equation}
which can be described by the commutative diagram
\begin{equation}
\begin{diagram} \left\{ \bullet \right\}\times G\times \left\{ \bullet \right\} & \rTo^{\mu_0^{\left(e\right)}\times\operatorname*{id}\times\mu_0^{\left(e\right)}} & G^{\times 3} & \rTo^{\operatorname*{id}^{\times 2}\times\bar{\mu}_{1}} & G^{\times 3} \\ \uTo^{\epsilon} & & & & \dTo_{\mu_{3}} \\ G & & \rTo^{ \varphi} & & G \\ \end{diagram}\label{dia-phi}%
\end{equation}
\smallskip

The mapping $\psi_{4}$ is the first polyadic power (\ref{pp}) of the binary
identity $e$ in the ternary system%
\begin{equation}
\psi_{4}\left(  e\right)  =e^{\left\langle 1\right\rangle }.\label{e4}%
\end{equation}
Thus, combining (\ref{m3h})--(\ref{e4}) we obtain the Hossz\'{u}-Gluskin
\textquotedblleft chain formula\textquotedblright\ for $n=3$%
\begin{align}
\mu_{3}\left[  g,h,u\right]    & =g\ast\varphi\left(  h\right)  \ast
\varphi^{2}\left(  u\right)  \ast b,\label{m3gh}\\
b  & =e^{\left\langle 1\right\rangle },\label{be11}%
\end{align}
which depends on one mapping $\varphi$ (taken in the chain of powers) only, and
the first polyadic power $e^{\left\langle 1\right\rangle }$ of the binary
identity $e$. The corresponding Hossz\'{u}-Gluskin diagram
\begin{equation}
\begin{diagram} G^{\times 3}\times \left\{ \bullet \right\}^3 & \rTo^{\operatorname*{id}\times\varphi\times\varphi^2\times\left(\mu_0^{\left(e\right)}\right)^{\times 3}} & G^{\times 6} & \rTo^{\operatorname*{id}^{\times 3}\times\mu_3} & G^{\times 4} \\ \uTo^{\epsilon} & & & & \dTo_{\mu_{2}^{\times 3}} \\ G\times G \times G & & \rTo^{ \mu_{3}} & & G \\ \end{diagram}\label{diaGH}%
\end{equation}
\noindent commutes.

The mapping $\varphi$ is an automorphism of the binary group $\left\langle
G\mid\ast,e\right\rangle $, because it follows from (\ref{g3}) and (\ref{fg})
that%
\begin{align}
\varphi\left(  g\right)  \ast\varphi\left(  h\right)   &  =\mu_{3}\left[
\mu_{3}\left[  e,g,\bar{e}\right]  ,\bar{e},\mu_{3}\left[  e,h,\bar{e}\right]
\right]  =\mathbf{\mu}_{3}^{3}\left[  e,g,\bar{e},\overset{\text{neutral}%
}{\left(  \bar{e},e\right)  },h,\bar{e}\right] \nonumber\\
&  =\mathbf{\mu}_{3}^{2}\left[  e,g,\bar{e},h,\bar{e}\right]  =\mu_{3}\left[
e,g\ast h,\bar{e}\right]  =\varphi\left(  g\ast h\right)  ,\\
\varphi\left(  e\right)   &  =\mu_{3}\left[  e,e,\bar{e}\right]  =\mu
_{3}\left[  e,\overset{\text{neutral}}{\left(  e,\bar{e}\right)  }\right]  =e.
\end{align}
It is important to note that not only the binary identity $e$, but also its
first polyadic power $e^{\left\langle 1\right\rangle }$ is a fixed point of
the automorphism $\varphi$, because%
\begin{equation}
\varphi\left(  e^{\left\langle 1\right\rangle }\right)  =\mu_{3}\left[
e,e^{\left\langle 1\right\rangle },\bar{e}\right]  =\mathbf{\mu}_{3}%
^{2}\left[  e,e,e,\overset{\text{neutral}}{\left(  e,\bar{e}\right)  }\right]
=\mu_{3}\left[  e,e,e\right]  =e^{\left\langle 1\right\rangle }. \label{fe1}%
\end{equation}
Moreover, taking into account that in the binary group (see (\ref{g3}))%
\begin{equation}
\left(  e^{\left\langle 1\right\rangle }\right)  ^{-1}=\mu_{3}\left[
e,\overline{e^{\left\langle 1\right\rangle }},e\right]  =\mathbf{\mu}_{3}%
^{2}\left[  e,\bar{e},\bar{e},\bar{e},e\right]  =\bar{e},
\end{equation}
we get
\begin{equation}
\varphi^{2}\left(  g\right)  =\mathbf{\mu}_{3}^{2}\left[  e,e,g,\bar{e}%
,\bar{e}\right]  =\mathbf{\mu}_{3}^{2}\left[  e,e,\overset{\text{neutral}%
}{\left(  e,\bar{e}\right)  }g,\bar{e},\bar{e}\right]  =e^{\left\langle
1\right\rangle }\ast g\ast\left(  e^{\left\langle 1\right\rangle }\right)
^{-1}. \label{f2}%
\end{equation}
The higher polyadic powers $e^{\left\langle k\right\rangle }=\mathbf{\mu}%
_{3}^{k}\left[  e^{2k+1}\right]  $ of the binary identity $e$ are obviously
also fixed points%
\begin{equation}
\varphi\left(  e^{\left\langle k\right\rangle }\right)  =e^{\left\langle
k\right\rangle }. \label{fek}%
\end{equation}

The elements $e^{\left\langle k\right\rangle }$ form a subgroup $\mathfrak{H}$
of the binary group $\left\langle G\mid\ast,e\right\rangle $, because%
\begin{align}
e^{\left\langle k\right\rangle }\ast e^{\left\langle l\right\rangle } &
=e^{\left\langle k+l\right\rangle },\\
e^{\left\langle k\right\rangle }\ast e &  =e\ast e^{\left\langle
k\right\rangle }=e^{\left\langle k\right\rangle }.
\end{align}
We can express the even powers of the automorphism $\varphi$ through the
polyadic powers $e^{\left\langle k\right\rangle }$ in the following way%
\begin{equation}
\varphi^{2k}\left(  g\right)  =e^{\left\langle k\right\rangle }\ast
g\ast\left(  e^{\left\langle k\right\rangle }\right)  ^{-1}.\label{fk}%
\end{equation}
This gives a manifest connection between the Hossz\'{u}-Gluskin
\textquotedblleft chain formula\textquotedblright\ and the sequence of cosets
(see, \cite{pos}) for the particular case $n=3$.

\begin{example}
Let us consider the ternary copula associative multiplication
\cite{mes/sar,stu/kol}%
\begin{equation}
\mu_{3}\left[  g,h,u\right]  =\dfrac{g(1-h)u}{g(1-h)u+(1-g)h(1-u)},
\end{equation}
where $g_{i}\in G=\left[  0,1\right]  $ and $0/0=0$ is assumed\footnote{In
this example all denominators are supposed nonzero.}. It is associative and
cannot be iterated from any binary group. Obviously,  $\mu_{3}\left[
g^{3}\right]  =g$, and therefore this polyadic system is $\ell_{\mu}%
$-idempotent (\ref{mi}) $g^{\left\langle \ell_{\mu}\right\rangle }=g$. The
querelement is $\bar{g}=\bar{\mu}_{1}\left[  g\right]  =g$. Because each
element is querable, then $\left\langle G\mid\mu_{3},\bar{\mu}_{1}%
\right\rangle $ is a ternary group. Take a fixed element $e\in\left[
0,1\right]  $. We define the binary multiplication as $g\ast h=\mu_{3}\left[
g,e,h\right]  $ and the automorphism
\begin{equation}
\varphi\left(  g\right)  =\mu_{3}\left[  e,g,e\right]  =e^{2}\dfrac{1-g}%
{e^{2}-2ge+g} \label{ege}%
\end{equation}
which has the property $\varphi^{2k}=\operatorname*{id}$ and $\varphi
^{2k+1}=\varphi$, where $k\in\mathbb{N}$. Obviously,  in (\ref{ege}) $g$
can be on any place in the product $\mu_{3}\left[  e,g,e\right]  =\mu
_{3}\left[  e,e,g\right]  =\mu_{3}\left[  e,e,g\right]  $. Now we can check
the Hossz\'{u}-Gluskin \textquotedblleft chain formula\textquotedblright%
\ (\ref{m3gh}) for the ternary copula%
\begin{align}
\mu_{3}\left[  g,h,u\right]   &  =\left(  \left(  \left(  g\ast\varphi\left(
h\right)  \right)  \ast u\right)  \ast e\right)  =\mathbf{\mu}_{3}^{\bullet
}\left[  g,e,e^{2}\dfrac{1-h}{e^{2}-2he+g},e,\left(  u,e,e\right)  \right]
\nonumber\\
&  =\mathbf{\mu}_{3}^{\bullet}\left[  g,\left(  e,e^{2}\dfrac{1-h}%
{e^{2}-2he+g},e\right)  ,u\right]  =\mu_{3}\left[  g,\varphi^{2}\left(
h\right)  ,u\right]  =\mu_{3}\left[  g,h,u\right]  .
\end{align}

\end{example}

The language of polyadic inverses allows us to generalize the
Hossz\'{u}-Gluskin \textquotedblleft chain formula\textquotedblright\ from
$n=3$ (\ref{m3gh}) to arbitrary $n$ in a clear way. The derivation coincides
with (\ref{m3h}) using the multiplication (\ref{ghe}) (with substitution
$\bar{e}\rightarrow\mathbf{e}^{-1}$), neutral polyads $\left(  \mathbf{e}%
^{-1},e\right)  $ or their powers $\left(  \left(  \mathbf{e}^{-1}\right)
^{k},e^{k}\right)  $, but contains $n$ terms%
\begin{align}
&  \mu_{n}\left[  g_{1},\ldots,g_{n}\right]  =\mathbf{\mu}_{n}^{\bullet
}\left[  g_{1},\overset{%
\begin{array}
[c]{c}%
\ast\\
\downarrow
\end{array}
}{\mathbf{e}^{-1}},e,g_{2},\ldots,g_{n}\right]  =\mathbf{\mu}_{n}^{\bullet
}\left[  g_{1},\overset{%
\begin{array}
[c]{c}%
\ast\\
\downarrow
\end{array}
}{\mathbf{e}^{-1}},e,g_{2},\mathbf{e}^{-1},\overset{%
\begin{array}
[c]{c}%
\ast\\
\downarrow
\end{array}
}{\mathbf{e}^{-1}},e,e,g_{3},\ldots,g_{n}\right]  =\ldots\nonumber\\
&  =\mathbf{\mu}_{n}^{\bullet}\left[  g_{1},\overset{%
\begin{array}
[c]{c}%
\ast\\
\downarrow
\end{array}
}{\mathbf{e}^{-1}},e,g_{2},\mathbf{e}^{-1},\overset{%
\begin{array}
[c]{c}%
\ast\\
\downarrow
\end{array}
}{\mathbf{e}^{-1}},e,e,g_{3},\ldots,\overset{%
\begin{array}
[c]{c}%
\ast\\
\downarrow
\end{array}
}{\mathbf{e}^{-1}},\overset{n-1}{\overbrace{e,\ldots,e}},g_{n},\overset
{n-1}{\overbrace{\mathbf{e}^{-1},\ldots,\mathbf{e}^{-1}}},\overset{%
\begin{array}
[c]{c}%
\ast\\
\downarrow
\end{array}
}{\mathbf{e}^{-1}},\overset{n}{\overbrace{e,\ldots,e}}\right]  . \label{mn1}%
\end{align}

We observe from (\ref{mn1}) that the mapping $\varphi$ in the $n$-ary case is%
\begin{equation}
\varphi\left(  g\right)  =\mu_{n}\left[  e,g,\mathbf{e}^{-1}\right]  ,
\label{fgn}%
\end{equation}
and the last product of the binary identities $\mu_{n}\left[  e,\ldots
,e\right]  $ is also the first $n$-ary power $e^{\left\langle 1\right\rangle
}$ (\ref{pp}). It follows from (\ref{fgn}) and (\ref{ghe}), that
\begin{equation}
\varphi^{n-1}\left(  g\right)  =e^{\left\langle 1\right\rangle }\ast
g\ast\left(  e^{\left\langle 1\right\rangle }\right)  ^{-1}. \label{fn1}%
\end{equation}
In this way, we obtain the Hossz\'{u}-Gluskin \textquotedblleft chain
formula\textquotedblright\ for arbitrary $n$%
\begin{equation}
\mu_{n}\left[  g_{1},\ldots,g_{n}\right]  =g_{1}\ast\varphi\left(
g_{2}\right)  \ast\varphi^{2}\left(  g_{3}\right)  \ast\ldots\ast\varphi
^{n-2}\left(  g_{n-1}\right)  \ast\varphi^{n-1}\left(  g_{n}\right)  \ast
e^{\left\langle 1\right\rangle }=\left(  \ast\underset{i=1}{\overset{n}{%
{\displaystyle\prod}
}}\varphi^{i-1}\left(  g_{i}\right)  \right)  \ast e^{\left\langle
1\right\rangle }. \label{gh}%
\end{equation}
Thus, we have found the \textquotedblleft extended\textquotedblright\ homotopy
maps $\psi_{i}$ from (\ref{mge}) as%
\begin{align}
\psi_{i}\left(  g\right)   &  =\varphi^{i-1}\left(  g\right)
,\ \ \ i=1,\ldots,n,\label{pg1}\\
\psi_{n+1}\left(  g\right)   &  =g^{\left\langle 1\right\rangle }, \label{pg2}%
\end{align}
where we put by definition $\varphi^{0}\left(  g\right)  =g$. Using
(\ref{fe1}) and (\ref{gh}) we can formulate the  Hossz\'{u}-Gluskin
theorem in the language of polyadic powers.

\begin{theorem}
\label{th-gh}On a polyadic group $\mathfrak{G}_{n}=\left\langle G\mid\mu
_{n},\bar{\mu}_{1}\right\rangle $ one can define a binary group $\mathfrak{G}%
_{2}^{\ast}=\left\langle G\mid\mu_{2}=\ast,e\right\rangle $ and its
automorphism $\varphi$ such that the Hossz\'{u}-Gluskin \textquotedblleft
chain formula\textquotedblright\ \textrm{(\ref{gh})} is valid, where the
polyadic powers of the identity $e$ are fixed points of $\varphi$
\textrm{(\ref{fek})}, form a subgroup $\mathfrak{H}$ of $\mathfrak{G}%
_{2}^{\ast}$, and the $\left(  n-1\right)  $ power of $\varphi$ is a
conjugation \textrm{(\ref{fn1})} with respect to $\mathfrak{H}$.
\end{theorem}

The following reverse Hossz\'{u}-Gluskin theorem holds.

\begin{theorem}
\label{th-ghr}If in a binary group $\mathfrak{G}_{2}^{\ast}=\left\langle
G\mid\mu_{2}=\ast,e\right\rangle $ one can define an automorphism $\varphi$
such that%
\begin{align}
\varphi^{n-1}\left(  g\right)   &  =b\ast g\ast b^{-1},\label{b1}\\
\varphi\left(  b\right)   &  =b,\label{b2}%
\end{align}
where $b\in G$ is a distinguished element, then the \textquotedblleft chain
formula\textquotedblright\
\begin{equation}
\mu_{n}\left[  g_{1},\ldots,g_{n}\right]  =\left(  \ast\underset{i=1}%
{\overset{n}{%
{\displaystyle\prod}
}}\varphi^{i-1}\left(  g_{i}\right)  \right)  \ast b\label{ghr}%
\end{equation}
determines a $n$-ary group, in which the distinguished element is the first
polyadic power of the binary identity%
\begin{equation}
b=e^{\left\langle 1\right\rangle }.\label{be1}%
\end{equation}

\end{theorem}

\section{\textquotedblleft Deformation\textquotedblright\ of
Hossz\'{u}-Gluskin chain formula}

Let us raise the question: can the choice (\ref{pg1})-(\ref{pg2}) of the
\textquotedblleft extended\textquotedblright\ homotopy maps (\ref{mge}) be
generalized? Before answering this question \textit{positively} we consider
some preliminary statements.

First, we note that we keep the general idea of inserting neutral sequences
into a polyadic product (see (\ref{m3g}) and (\ref{mn1})), because this is the
only way to obtain \textquotedblleft automatic\textquotedblright%
\ associativity. Second, the number of the inserted neutral polyads can be
chosen \textit{arbitrarily}, not only minimally, as in (\ref{m3g}) and
(\ref{mn1}) (as they are neutral). Nevertheless, we can show that this
arbitrariness is somewhat restricted.

Indeed, let us consider a polyadic group $\left\langle G\mid\mu_{n},\bar{\mu
}_{1}\right\rangle $ in the particular case $n=3$, where for any $e_{0}\in G$
and natural $k$ the sequence $\left(  \bar{e}_{0}^{k},e_{0}^{k}\right)  $ is
neutral, then we can write%
\begin{equation}
\mu_{3}\left[  g,h,u\right]  =\mathbf{\mu}_{3}^{\bullet}\left[  g,\bar{e}%
_{0}^{k},e_{0}^{k},h,\bar{e}_{0}^{lk},e_{0}^{lk},u,\bar{e}_{0}^{mk},e_{0}%
^{mk}\right]  .
\end{equation}
If we make the change of variables $e_{0}^{k}=e$, then we obtain%
\begin{equation}
\mu_{3}\left[  g,h,u\right]  =\mathbf{\mu}_{3}^{\bullet}\left[  g,\bar
{e},e,h,\bar{e}^{l},e^{l},u,\bar{e}^{m},e^{m}\right]  . \label{me1}%
\end{equation}
Because this should reproduce the formula (\ref{m3}), we immediately conclude
that $\psi_{1}\left(  g\right)  =\operatorname*{id}$, and the multiplication
is the same as in (\ref{g3}), and $e$ is again the identity of the binary
group $\mathfrak{G}^{\ast}=\left\langle G,\ast,e\right\rangle $. Moreover, if
we put $\psi_{2}\left(  g\right)  =\varphi\left(  g\right)  $, as in the
standard case, then we have a first \textquotedblleft half\textquotedblright%
\ of the mapping $\varphi$, that is $\varphi\left(  g\right)  =\mu_{3}\left[
e,h,\text{something}\right]  $. Now we are in a position to find this
\textquotedblleft something\textquotedblright\ and other \textquotedblleft
extended\textquotedblright\ homotopy maps $\psi_{i}$ from (\ref{m3}), but
\textit{without} the requirement of a minimal number of inserted neutral
polyads, as it was in (\ref{m3g}). By analogy, we rewrite (\ref{me1}) as%
\begin{equation}
\mu_{3}\left[  g,h,u\right]  =\mathbf{\mu}_{3}^{\bullet}\left[  g,\bar
{e},\left(  e,h,\bar{e}^{q}\right)  ,\bar{e},e^{q+1},u,\bar{e}^{m}%
,e^{m}\right]  , \label{mgq}%
\end{equation}
where we put $l=q+1$. So we have found the \textquotedblleft
something\textquotedblright, and the map $\varphi$ is%
\begin{equation}
\varphi_{q}\left(  g\right)  =\mathbf{\mu}_{3}^{\ell_{\varphi}\left(
q\right)  }\left[  e,g,\bar{e}^{q}\right]  , \label{fgq}%
\end{equation}
where the number of multiplications%
\begin{equation}
\ell_{\varphi}\left(  q\right)  =\dfrac{q+1}{2} \label{kf}%
\end{equation}
is an integer $\ell_{\varphi}\left(  q\right)  =1,2,3\ldots$, while
$q=1,3,5,7\ldots$. The diagram defined $\varphi_{q}$ (e.g., for $q=3$ and
$\ell_{\varphi}\left(  q\right)  =2$)
\begin{equation}
\begin{diagram} \left\{ \bullet \right\}\times G\times \left\{ \bullet \right\}^3 & \rTo^{\mu_0^{\left(e\right)}\times\operatorname*{id}\times\left(\mu_0^{\left(e\right)}\right)^3} & G^{\times 5} & \rTo^{\operatorname*{id}^{\times 2}\times\left(\bar{\mu}_{1}\right)^3} & G^{\times 5} \\ \uTo^{\epsilon} & & & & \dTo_{\mu_{3}\times\mu_{3}} \\ G & & \rTo^{ \varphi_q} & & G \\ \end{diagram}\label{dia-phiq}%
\end{equation}
commutes (cf. (\ref{dia-phi})). Then, we can find power $m$ in (\ref{mgq})%
\begin{equation}
\mu_{3}\left[  g,h,u\right]  =\mathbf{\mu}_{3}^{\bullet}\left[  g,\bar
{e},\left(  e,h,\bar{e}^{q}\right)  ,\bar{e},\left(  e,u,\bar{e}^{q}\right)
^{q+1},\bar{e},e^{q\left(  q+1\right)  +1}\right]  ,
\end{equation}
and therefore $m=q\left(  q+1\right)  +1$. Thus, we have obtained the
\textquotedblleft$q$-deformed\textquotedblright\ maps $\psi_{i}$ (cf.
(\ref{p1})--(\ref{p4}))%

\begin{align}
\psi_{1}\left(  g\right)   &  =\varphi_{q}^{\left[  \left[  0\right]  \right]
_{q}}\left(  g\right)  =\varphi_{q}^{0}\left(  g\right)  =g,\label{f1}\\
\psi_{2}\left(  g\right)   &  =\varphi_{q}\left(  g\right)  =\varphi
_{q}^{\left[  \left[  1\right]  \right]  _{q}}\left(  g\right)  ,\\
\psi_{3}\left(  g\right)   &  =\varphi_{q}^{q+1}\left(  g\right)  =\varphi
_{q}^{\left[  \left[  2\right]  \right]  _{q}}\left(  g\right)  ,\\
\psi_{4}\left(  g\right)   &  =\mathbf{\mu}_{3}^{\bullet}\left[  g^{q\left(
q+1\right)  +1}\right]  =\mathbf{\mu}_{3}^{\bullet}\left[  g^{\left[  \left[
3\right]  \right]  _{q}}\right]  ,\label{f4}%
\end{align}
where $\varphi$ is defined by (\ref{fgq}) and $\left[  \left[  k\right]
\right]  _{q}$ is the $q$-deformed number (\ref{kq}), and we put $\varphi
_{q}^{0}=\operatorname*{id}$. The corresponding \textquotedblleft%
$q$-deformed\textquotedblright\ chain formula (for $n=3$) can be written as
(cf. (\ref{m3gh})--(\ref{be11}) for \textquotedblleft
nondeformed\textquotedblright\ case)%
\begin{align}
\mu_{3}\left[  g,h,u\right]    & =g\ast\varphi_{q}^{\left[  \left[  1\right]
\right]  _{q}}\left(  h\right)  \ast\varphi_{q}^{\left[  \left[  2\right]
\right]  _{q}}\left(  u\right)  \ast b_{q},\label{m3q}\\
b_{q}  & =e^{\left\langle \ell_{e}\left(  q\right)  \right\rangle },\label{be}%
\end{align}
where the degree of the binary identity polyadic power%
\begin{equation}
\ell_{e}\left(  q\right)  =q\dfrac{\left[  \left[  2\right]  \right]  _{q}}%
{2}=\ell_{\varphi}\left(  q\right)  \left(  2\ell_{\varphi}\left(  q\right)
+1\right)
\end{equation}
is an integer. The corresponding \textquotedblleft deformed\textquotedblright%
\ chain diagram (e.g., for $q=3$)
\begin{equation}
\begin{diagram} G^{\times 3}\times \left\{ \bullet \right\}^{13} & \rTo^{\operatorname*{id}\times\varphi_q\times\varphi_q^4\times\left(\mu_0^{\left(e\right)}\right)^{\times 13}} & G^{\times 16} & \rTo^{\operatorname*{id}^{\times 3}\times\mu_3^6} & G^{\times 4} \\ \uTo^{\epsilon} & & & & \dTo_{\mu_{2}^{\times 3}} \\ G\times G \times G & & \rTo^{ \mu_{3}} & & G \\ \end{diagram}
\end{equation}
\noindent commutes (cf. the Hossz\'{u}-Gluskin diagram (\ref{diaGH})). In the
\textquotedblleft deformed\textquotedblright\ case the polyadic power
$e^{\left\langle \ell_{e}\left(  q\right)  \right\rangle }$ is not a fixed
point of $\varphi_{q}$ and satisfies%
\begin{equation}
\varphi_{q}\left(  e^{\left\langle \ell_{e}\left(  q\right)  \right\rangle
}\right)  =\varphi_{q}\left(  \mathbf{\mu}_{3}^{\bullet}\left[  e^{q^{2}%
+q+1}\right]  \right)  =\mathbf{\mu}_{3}^{\bullet}\left[  e^{q^{2}+2}\right]
=e^{\left\langle \ell_{e}\left(  q\right)  \right\rangle }\ast\varphi
_{q}\left(  e\right)  \label{fq1}%
\end{equation}
or%
\begin{equation}
\varphi_{q}\left(  b_{q}\right)  =b_{q}\ast\varphi_{q}\left(  e\right)  .
\end{equation}
Instead of (\ref{f2}) we have%
\begin{equation}
\varphi_{q}^{q+1}\left(  g\right)  \ast e^{\left\langle \ell_{e}\left(
q\right)  \right\rangle }=\mathbf{\mu}_{3}^{\bullet}\left[  e^{q+1},g\right]
=\mathbf{\mu}_{3}^{\bullet}\left[  e^{q+2}\right]  \ast g=e^{\left\langle
\ell_{e}\left(  q\right)  \right\rangle }\ast\varphi_{q}^{q+1}\left(
e\right)  \ast g\label{fqq}%
\end{equation}
or%
\begin{equation}
\varphi_{q}^{q+1}\left(  g\right)  \ast b_{q}=b_{q}\ast\varphi_{q}%
^{q+1}\left(  e\right)  \ast g.
\end{equation}

The \textquotedblleft nondeformed\textquotedblright\ limit $q\rightarrow1$ of
(\ref{m3q}) gives the Hossz\'{u}-Gluskin chain formula (\ref{m3gh})
for $n=3$. Now let us turn to arbitrary $n$ and write the $n$-ary
multiplication using neutral polyads analogously to (\ref{mgq}). By the same
arguments, as in (\ref{me1}), we insert only one neutral polyad $\left(
\mathbf{e}^{-1},e\right)  $ between the first and second elements in the
multiplication, but in other places we insert powers $\left(  \left(
\mathbf{e}^{-1}\right)  ^{k},e^{k}\right)  $ (\textit{allowed} by the chain
properties), and obtain%
\begin{align}
&  \mu_{n}\left[  g_{1},\ldots,g_{n}\right]  =\mathbf{\mu}_{n}^{\bullet
}\left[  g_{1},\mathbf{e}^{-1},e,g_{2},\ldots,g_{n}\right]  =\mathbf{\mu}%
_{n}^{\bullet}\left[  g_{1},\mathbf{e}^{-1},\left(  e,g_{2},\left(
\mathbf{e}^{-1}\right)  ^{q}\right)  ,\mathbf{e}^{-1},e^{q+1},g_{3}%
,\ldots,g_{n}\right]  =\ldots\nonumber\\
&  =\mathbf{\mu}_{n}^{\bullet}\left[  g_{1},\mathbf{e}^{-1},\left(
e,g_{2},\left(  \mathbf{e}^{-1}\right)  ^{q}\right)  ,\mathbf{e}^{-1},\left(
e^{q+1},g_{3},\overset{q\left(  q+1\right)  }{\overbrace{\mathbf{e}%
^{-1},\ldots,\mathbf{e}^{-1}}}\right)  \mathbf{e}^{-1},e^{q\left(  q+1\right)
+1},g_{3},\ldots\right. \nonumber\\
&  \left.  \ldots,\left(  \overset{q^{n-2}+\ldots+q+1}{\overbrace{e,\ldots,e}%
},g_{n-1},\overset{q\left(  q^{n-2}+\ldots+q+1\right)  }{\overbrace
{\mathbf{e}^{-1},\ldots,\mathbf{e}^{-1}}}\right)  ,\mathbf{e}^{-1},\left(
\overset{q^{n-1}+\ldots+q+1}{\overbrace{e,\ldots,e}},g_{n},\overset{q\left(
q^{n-1}+\ldots+q+1\right)  }{\overbrace{\mathbf{e}^{-1},\ldots,\mathbf{e}%
^{-1}}}\right)  ,\mathbf{e}^{-1},\overset{q^{n}+\ldots+q+1}{\overbrace
{e,\ldots,e}}\right]  . \label{mnn}%
\end{align}
So we observe that the binary product is now the same as in the
\textquotedblleft nondeformed\textquotedblright\ case (\ref{ghe}), while the
map $\varphi$ is%
\begin{equation}
\varphi_{q}\left(  g\right)  =\mathbf{\mu}_{n}^{\ell_{\varphi}\left(
q\right)  }\left[  e,g,\left(  \mathbf{e}^{-1}\right)  ^{q}\right]  ,
\label{fq}%
\end{equation}
where the number of multiplications%
\begin{equation}
\ell_{\varphi}\left(  q\right)  =\dfrac{q\left(  n-2\right)  +1}{n-1}
\label{lf}%
\end{equation}
is an integer and $\ell_{\varphi}\left(  q\right)  \rightarrow q$, as
$n\rightarrow\infty$, in the nondeformed case $\ell_{\varphi}\left(  1\right)
=1$, as in (\ref{fgn}). Note that the \textquotedblleft
deformed\textquotedblright\ map $\varphi_{q}$ is the $a$-quasi-endomorphism
\cite{glu/shv} of the binary group $\mathfrak{G}_{2}^{\ast}$, because from
(\ref{fq}) we get%
\begin{align}
\varphi_{q}\left(  g\right)  \ast\varphi_{q}\left(  h\right)   &
=\mathbf{\mu}_{n}^{\bullet}\left[  e,g,\left(  \mathbf{e}^{-1}\right)
^{q},\mathbf{e}^{-1},e,h,\left(  \mathbf{e}^{-1}\right)  ^{q}\right]
\nonumber\\
&  =\mathbf{\mu}_{n}^{\bullet}\left[  e,g,\mathbf{e}^{-1},\left(  e,e,\left(
\mathbf{e}^{-1}\right)  ^{q}\right)  ,\mathbf{e}^{-1},h,\left(  \mathbf{e}%
^{-1}\right)  ^{q}\right]  =\varphi_{q}\left(  g\ast a\ast h\right)  ,
\end{align}
where%
\begin{equation}
a=\mathbf{\mu}_{n}^{\ell_{\varphi}\left(  q\right)  }\left[  e,e,\left(
\mathbf{e}^{-1}\right)  ^{q}\right]  =\varphi_{q}\left(  e\right)  . \label{a}%
\end{equation}
In general, a \textit{quasi-endomorphism} can be defined by%
\begin{equation}
\varphi_{q}\left(  g\right)  \ast\varphi_{q}\left(  h\right)  =\varphi
_{q}\left(  g\ast\varphi_{q}\left(  e\right)  \ast h\right)  . \label{ff}%
\end{equation}
The corresponding diagram
\begin{equation}
\begin{diagram} G\times G & \rTo^{\mu_2} & G & \lTo^{\varphi_q} & G \\ \uTo^{\varphi_q\times\varphi_q} & & & & \uTo_{\mu_{2}\times \mu_2} \\ G\times G &\rTo^{\epsilon} &G\times \left\{ \bullet \right\}\times G&\rTo^{ \operatorname*{id}\times\mu_0^{\left(e\right)}\times\operatorname*{id}} & G \times G\times G\\ \end{diagram}\label{dia-quasi}%
\end{equation}
\noindent commutes. If $q=1$, then $\varphi_{q}\left(  e\right)  =e$, and the
distinguished element $a$ turns to the binary identity $a=e$, such that the
$a$-quasi-endomorphism $\varphi_{q}$ becomes an automorphism of $\mathfrak{G}%
_{2}^{\ast}$.

\begin{remark}
The choice (\ref{fq}) of the $a$-quasi-endomorphism $\varphi_{q}$ is different
from \cite{glu/shv}, the latter (in our notation) is $\varphi_{k}\left(
g\right)  =\mu_{n}\left[  a^{k-1},g,a^{n-k}\right]  $, $k=1,\ldots,n-1$, it
has only one multiplication and leads to the  \textquotedblleft
nondeformed\textquotedblright\ chain formula (\ref{gh}) (for semigroup case).
\end{remark}

It follows from (\ref{mnn}), that the \textquotedblleft
extended\textquotedblright\ homotopy maps $\psi_{i}$ (\ref{mge}) are (cf.
(\ref{f1})--(\ref{f4}))%
\begin{align}
\psi_{1}\left(  g\right)   &  =\varphi_{q}^{\left[  \left[  0\right]  \right]
_{q}}\left(  g\right)  =\varphi_{q}^{0}\left(  g\right)  =g,\label{pq}\\
\psi_{2}\left(  g\right)   &  =\varphi_{q}\left(  g\right)  =\varphi
_{q}^{\left[  \left[  1\right]  \right]  _{q}}\left(  g\right)  ,\\
\psi_{3}\left(  g\right)   &  =\varphi_{q}^{q+1}\left(  g\right)  =\varphi
_{q}^{\left[  \left[  2\right]  \right]  _{q}}\left(  g\right)  ,\\
&  \vdots\nonumber\\
\psi_{n-1}\left(  g\right)   &  =\varphi_{q}^{q^{n-3}+\ldots+q+1}\left(
g\right)  =\varphi_{q}^{\left[  \left[  n-2\right]  \right]  _{q}}\left(
g\right)  ,\\
\psi_{n}\left(  g\right)   &  =\varphi_{q}^{q^{n-2}+\ldots+q+1}\left(
g\right)  =\varphi_{q}^{\left[  \left[  n-1\right]  \right]  _{q}}\left(
g\right)  ,\\
\psi_{n+1}\left(  g\right)   &  =\mathbf{\mu}_{n}^{\bullet}\left[
g^{q^{n-1}+\ldots+q+1}\right]  =\mathbf{\mu}_{n}^{\bullet}\left[  g^{\left[
\left[  n\right]  \right]  _{q}}\right]  . \label{pq1}%
\end{align}
In terms of the polyadic power (\ref{pp}), the last map is%
\begin{equation}
\psi_{n+1}\left(  g\right)  =g^{\left\langle \ell_{e}\right\rangle },
\end{equation}
where (cf. (\ref{lf}))%
\begin{equation}
\ell_{e}\left(  q\right)  =q\dfrac{\left[  \left[  n-1\right]  \right]  _{q}%
}{n-1} \label{leq}%
\end{equation}
is an integer. Thus the \textquotedblleft$q$-deformed\textquotedblright%
\ $n$-ary chain formula is (cf. (\ref{gh}))%
\begin{equation}
\mu_{n}\left[  g_{1},\ldots,g_{n}\right]  =g_{1}\ast\varphi_{q}^{\left[
\left[  1\right]  \right]  _{q}}\left(  g_{2}\right)  \ast\varphi_{q}^{\left[
\left[  2\right]  \right]  _{q}}\left(  g_{3}\right)  \ast\ldots\ast
\varphi_{q}^{\left[  \left[  n-2\right]  \right]  _{q}}\left(  g_{n-1}\right)
\ast\varphi_{q}^{\left[  \left[  n-1\right]  \right]  _{q}}\left(
g_{n}\right)  \ast e^{\left\langle \ell_{e}\left(  q\right)  \right\rangle }.
\label{ghq}%
\end{equation}

In the \textquotedblleft nondeformed\textquotedblright\ limit $q\rightarrow1$
(\ref{ghq}) reproduces the  Hossz\'{u}-Gluskin chain formula
(\ref{gh}). Let us obtain the \textquotedblleft deformed\textquotedblright%
\ analogs of the distinguished element relations (\ref{b1})--(\ref{b2}) for
arbitrary $n$ (the case $n=3$ is in (\ref{fq1})--(\ref{fqq})). Instead of the
fixed point relation (\ref{b2}) we now have from (\ref{fq}), (\ref{leq}) and
(\ref{pq1}) the \textit{quasi-fixed point}%
\begin{equation}
\varphi_{q}\left(  b_{q}\right)  =b_{q}\ast\varphi_{q}\left(  e\right)  ,
\label{fbf}%
\end{equation}
where the \textquotedblleft deformed\textquotedblright\ distinguished element
$b_{q}$ is (cf. (\ref{be1}))%
\begin{equation}
b_{q}=\mathbf{\mu}_{n}^{\bullet}\left[  e^{\left[  \left[  n\right]  \right]
_{q}}\right]  =e^{\left\langle \ell_{e}\left(  q\right)  \right\rangle }.
\label{beq}%
\end{equation}

The conjugation relation (\ref{b1}) in the \textquotedblleft
deformed\textquotedblright\ case becomes the \textit{quasi-conjugation}%
\begin{equation}
\varphi_{q}^{\left[  \left[  n-1\right]  \right]  _{q}}\left(  g\right)  \ast
b_{q}=b_{q}\ast\varphi_{q}^{\left[  \left[  n-1\right]  \right]  _{q}}\left(
e\right)  \ast g.\label{fnq}%
\end{equation}
This allows us to rewrite the \textquotedblleft deformed\textquotedblright%
\ chain formula\ (\ref{ghq}) as%
\begin{equation}
\mu_{n}\left[  g_{1},\ldots,g_{n}\right]  =g_{1}\ast\varphi_{q}^{\left[
\left[  1\right]  \right]  _{q}}\left(  g_{2}\right)  \ast\varphi_{q}^{\left[
\left[  2\right]  \right]  _{q}}\left(  g_{3}\right)  \ast\ldots\ast
\varphi_{q}^{\left[  \left[  n-2\right]  \right]  _{q}}\left(  g_{n-1}\right)
\ast b_{q}\ast\varphi_{q}^{\left[  \left[  n-1\right]  \right]  _{q}}\left(
e\right)  \ast g_{n}.\label{mngg}%
\end{equation}
Using the above proof sketch, we formulate the following \textquotedblleft%
$q$-deformed\textquotedblright\ analog of the Hossz\'{u}-Gluskin theorem:

\begin{theorem}
On a polyadic group $\mathfrak{G}_{n}=\left\langle G\mid\mu_{n},\bar{\mu}%
_{1}\right\rangle $ one can define a binary group $\mathfrak{G}_{2}^{\ast
}=\left\langle G\mid\mu_{2}=\ast,e\right\rangle $ and (the infinite
\textquotedblleft$q$-series\textquotedblright\ of) its automorphism
$\varphi_{q}$ such that the \textquotedblleft deformed\textquotedblright%
\ chain formula\ \textrm{(\ref{ghq})} is valid%
\begin{equation}
\mu_{n}\left[  g_{1},\ldots,g_{n}\right]  =\left(  \ast\underset{i=1}%
{\overset{n}{%
{\displaystyle\prod}
}}\varphi^{\left[  \left[  i-1\right]  \right]  _{q}}\left(  g_{i}\right)
\right)  \ast b_{q},\label{mfq}%
\end{equation}
where (the infinite \textquotedblleft$q$-series\textquotedblright\ of) the
\textquotedblleft deformed\textquotedblright\ distinguished element $b_{q}$
(being a polyadic power of the binary identity \textrm{(\ref{beq})}) is the
quasi-fixed point of $\varphi_{q}$ \textrm{(\ref{fbf})} and satisfies the
quasi-conjugation \textrm{(\ref{fnq})} in the form%
\begin{equation}
\varphi_{q}^{\left[  \left[  n-1\right]  \right]  _{q}}\left(  g\right)
=b_{q}\ast\varphi_{q}^{\left[  \left[  n-1\right]  \right]  _{q}}\left(
e\right)  \ast g\ast b_{q}^{-1}.
\end{equation}

\end{theorem}

In the \textquotedblleft nondeformed\textquotedblright\ case $q=1$ we obtain
the  Hossz\'{u}-Gluskin chain formula (\ref{gh}) and the corresponding
\textbf{Theorem \ref{th-gh}}.

\begin{example}
Let us have a binary group $\left\langle G\mid\left(  \cdot\right)
,1\right\rangle $ and a distinguished element $e\in G$, $e\neq1$, then we can
define a binary group $\mathfrak{G}_{2}^{\ast}=\left\langle G\mid\left(
\ast\right)  ,e\right\rangle $ by the product%
\begin{equation}
g\ast h=g\cdot e^{-1}\cdot h.
\end{equation}
The quasi-endomorphism%
\begin{equation}
\varphi_{q}\left(  g\right)  =e\cdot g\cdot e^{-q}%
\end{equation}
satisfies (\ref{ff}) with $\varphi_{q}\left(  e\right)  =e^{2-q}$, and we take%
\begin{equation}
b_{q}=e^{\left[  \left[  n\right]  \right]  _{q}}.
\end{equation}
Then we can obtain the \textquotedblleft$q$-deformed\textquotedblright\ chain
formula (\ref{mfq}) (for $q=1$ see, e.g., \cite{gal/vor}).
\end{example}

We observe that the chain formula is the \textquotedblleft$q$%
-series\textquotedblright\ of equivalence relations (\ref{mfq}), which can be
formulated as an invariance. Indeed, let us denote the r.h.s. of (\ref{mfq}) by
$\mathcal{M}_{q}\left(  g_{1},\ldots,g_{n}\right)  $, and the l.h.s. as
$\mathcal{M}_{0}\left(  g_{1},\ldots,g_{n}\right)  $, then the chain formula
can be written as some invariance (cf. associativity as an invariance (\ref{ghu})).

\begin{theorem}
On a polyadic group $\mathfrak{G}_{n}=\left\langle G\mid\mu_{n},\bar{\mu}%
_{1}\right\rangle $ we can define a binary group $\mathfrak{G}^{\ast
}=\left\langle G\mid\mu_{2}=\ast,e\right\rangle $ such that the following
invariance is valid%
\begin{equation}
\mathcal{M}_{q}\left(  g_{1},\ldots,g_{n}\right)
=invariant,\ \ \ \ \ q=0,1,\ldots\ , \label{mq}%
\end{equation}
where%
\begin{equation}
\mathcal{M}_{q}\left(  g_{1},\ldots,g_{n}\right)  =\left\{
\begin{array}
[c]{c}%
\mu_{n}\left[  g_{1},\ldots,g_{n}\right]  ,\ \ \ q=0,\\
\left(  \ast\underset{i=1}{\overset{n}{%
{\displaystyle\prod}
}}\varphi^{\left[  \left[  i-1\right]  \right]  _{q}}\left(  g_{i}\right)
\right)  \ast b_{q},\ \ \ q>0,
\end{array}
\right.
\end{equation}
and the distinguished element $b_{q}\in G$ and the quasi-endomorphism
$\varphi_{q}$ of $\mathfrak{G}_{2}^{\ast}$ are defined in \textrm{(\ref{beq})}
and \textrm{(\ref{fq})} respectively.
\end{theorem}

\begin{example}
Let us consider the ternary $q$-product used in the nonextensive statistics
\cite{niv/lem/wan}%
\begin{equation}
\mu_{3}\left[  g,t,u\right]  =\left(  g^{\hbar}+t^{\hbar}+u^{\hbar}-3\right)
^{\tfrac{1}{\hbar}}, \label{gt}%
\end{equation}
where $\hbar=1-q_{0}$, and $g,t,u\in G=\mathbb{R}_{+}$, $0<q_{0}<1$, and also
$g^{\hbar}+t^{\hbar}+u^{\hbar}-3>0$ (as for other terms inside brackets with
power $\tfrac{1}{\hbar}$ below). In case $\hbar\rightarrow0$ the $q$-product
becomes an iterated product in $\mathbb{R}_{+}$ as $\mu_{3}\left[
g,t,u\right]  \rightarrow gtu$. The quermap $\bar{\mu}_{1}$ is given by%
\begin{equation}
\bar{g}=\left(  3-g^{\hbar}\right)  ^{\tfrac{1}{\hbar}}. \label{3g}%
\end{equation}
The polyadic system $\mathfrak{G}_{n}=\left\langle G\mid\mu_{3},\bar{\mu}%
_{1}\right\rangle $ is a ternary group, because each element is querable. Take
a distinguished element $e\in G$ and use (\ref{g3}), (\ref{gt}) and (\ref{3g})
to define the product%
\begin{equation}
g\ast t=\left(  g^{\hbar}-e^{\hbar}+t^{\hbar}\right)  ^{\tfrac{1}{\hbar}}%
\end{equation}
of the binary group $\mathfrak{G}_{2}^{\ast}=\left\langle G\mid\mu_{2}=\left(
\ast\right)  ,e\right\rangle $.

1) \textsl{The Hossz\'{u}-Gluskin chain formula} ($q=1$). The
automorphism (\ref{fg}) of $\mathfrak{G}^{\ast}$ is now the identity map
$\varphi=\operatorname*{id}$. The first polyadic power of the distinguished
element $e$ is%
\begin{equation}
b=e^{\left\langle 1\right\rangle }=\mu_{3}\left[  e^{3}\right]  =\left(
3e^{\hbar}-3\right)  ^{\tfrac{1}{\hbar}}.
\end{equation}
The chain formula (\ref{m3gh}) can be checked as follows%
\begin{align}
\mu_{3}\left[  g,t,u\right]   &  =\left(  \left(  \left(  g\ast t\right)  \ast
u\right)  \ast b\right)  =\left(  \left(  \left(  g^{\hbar}-e^{\hbar}%
+t^{\hbar}\right)  -e^{\hbar}+u^{\hbar}\right)  -e^{\hbar}+b^{\hbar}\right)
^{\tfrac{1}{\hbar}}\nonumber\\
&  =\left(  g^{\hbar}-e^{\hbar}+t^{\hbar}-e^{\hbar}+u^{\hbar}-e^{\hbar
}+3e^{\hbar}-3\right)  ^{\tfrac{1}{\hbar}}=\left(  g^{\hbar}+t^{\hbar
}+u^{\hbar}-3\right)  ^{\tfrac{1}{\hbar}}.
\end{align}
2) \textsl{The \textquotedblleft}$q$-\textsl{deformed\textquotedblright\ chain
formula} (for conciseness we consider only the case $q=3$). Now the
quasi-endomorphism $\varphi_{q}$ (\ref{fgq}) is not the identity, but is%
\begin{equation}
\varphi_{q=3}\left(  g\right)  =\left(  g^{\hbar}-2e^{\hbar}+3\right)
^{\tfrac{1}{\hbar}}.
\end{equation}
In case $q=3$ we need its $4$th ($=q+1$) power (\ref{m3q})%
\begin{equation}
\varphi_{q=3}^{4}\left(  g\right)  =\left(  g^{\hbar}-8e^{\hbar}+12\right)
^{\tfrac{1}{\hbar}}.
\end{equation}
The deformed polyadic power $e^{\left\langle \ell_{e}\right\rangle }$ from
(\ref{m3q}) is (see, also, (\ref{f4}))%
\begin{equation}
b_{q=3}=e^{\left\langle 5\right\rangle }=\mathbf{\mu}_{3}^{5}\left[
e^{13}\right]  =\left(  13e^{\hbar}-18\right)  ^{\tfrac{1}{\hbar}}.
\end{equation}
Now we check the \textquotedblleft$q$-deformed\textquotedblright\ chain
formula (\ref{m3q}) as%
\begin{align}
\mu_{3}\left[  g,t,u\right]   &  =g\ast\varphi_{q=3}\left(  t\right)
\ast\varphi_{q=3}^{4}\left(  u\right)  \ast b_{q=3}=\left(  \left(  \left(
g\ast\varphi_{q=3}\left(  t\right)  \right)  \ast\varphi_{q=3}^{4}\left(
u\right)  \right)  \ast b_{q=3}\right) \\
&  =\left(  g^{\hbar}-e^{\hbar}+\left(  t^{\hbar}-2e^{\hbar}+3\right)
-e^{\hbar}+\left(  u^{\hbar}-8e^{\hbar}+12\right)  -e^{\hbar}+\left(
13e^{\hbar}-18\right)  \right)  ^{\tfrac{1}{\hbar}}\\
&  =\left(  g^{\hbar}+t^{\hbar}+u^{\hbar}-3\right)  ^{\tfrac{1}{\hbar}}.
\end{align}

In a similar way, one can check the \textquotedblleft$q$%
-deformed\textquotedblright\ chain formula for any allowed $q$ (determined by
(\ref{lf}) and (\ref{leq}) to obtain an infinite $q$-series of the chain
representation of the same $n$-ary multiplication.
\end{example}

\section{Generalized \textquotedblleft deformed\textquotedblright\ version of
the homomorphism theorem}

Let us consider a homomorphism of the binary groups entering into the
\textquotedblleft deformed\textquotedblright\ chain formula (\ref{mfq}) as
$\Phi^{\ast}:\mathfrak{G}_{2}^{\ast}\rightarrow\mathfrak{G}_{2}^{\ast\prime}$,
where $\mathfrak{G}_{2}^{\ast\prime}=\left\langle G^{\prime}\mid\ast^{\prime
},e^{\prime}\right\rangle $. We observe that, because $\Phi^{\ast}$ commutes
with the binary multiplication, we need its commutation also with the
automorphisms $\varphi_{q}$ in each term of (\ref{mfq}) (which fixes equality
of the \textquotedblleft deformation\textquotedblright\ parameters
$q=q^{\prime}$) and its homomorphic action on $b_{q}$. Indeed, if
\begin{align}
\Phi^{\ast}\left(  \varphi_{q}\left(  g\right)  \right)   &  =\varphi
_{q}^{\prime}\left(  \Phi^{\ast}\left(  g\right)  \right)  ,\label{ff1}\\
\Phi^{\ast}\left(  b_{q}\right)   &  =b_{q}^{\prime},\label{ff2}%
\end{align}
then we get from (\ref{mfq})%
\begin{align}
\Phi^{\ast}\left(  \mu_{n}\left[  g_{1},\ldots,g_{n}\right]  \right)   &
=\Phi^{\ast}\left(  g_{1}\right)  \ast^{\prime}\Phi^{\ast}\left(  \varphi
_{q}^{\left[  \left[  1\right]  \right]  _{q}}\left(  g_{2}\right)  \right)
\ast^{\prime}\ldots\ast^{\prime}\Phi^{\ast}\left(  \varphi_{q}^{\left[
\left[  n-1\right]  \right]  _{q}}\left(  g_{n}\right)  \right)  \ast^{\prime
}\Phi^{\ast}\left(  b_{q}\right)  \nonumber\\
&  =\Phi^{\ast}\left(  g_{1}\right)  \ast^{\prime}\varphi_{q}^{\prime\left[
\left[  1\right]  \right]  _{q}}\left(  \Phi^{\ast}\left(  g_{2}\right)
\right)  \ast^{\prime}\ldots\ast^{\prime}\varphi_{q}^{\left[  \left[
n-1\right]  \right]  _{q}}\left(  \Phi^{\ast}\left(  g_{n}\right)  \right)
\ast^{\prime}b_{q}^{\prime}\nonumber\\
&  =\mu_{n}^{\prime}\left[  \Phi^{\ast}\left(  g_{1}\right)  ,\ldots
,\Phi^{\ast}\left(  g_{n}\right)  \right]  ,\label{ff3}%
\end{align}
where $g^{\prime}\ast^{\prime}h^{\prime}=\mu_{n}^{\prime}\left[  g^{\prime
},\mathbf{e}^{\prime-1},h^{\prime}\right]  $, $\varphi_{q}^{\prime}\left(
g^{\prime}\right)  =\mathbf{\mu}_{n}^{\prime\ \ell_{\varphi}\left(  q\right)
}\left[  e^{\prime},g^{\prime},\left(  \mathbf{e}^{\prime-1}\right)
^{q}\right]  $, $b_{q}^{\prime}=\mathbf{\mu}_{n}^{\prime\bullet}\left[
e^{\prime\ \left[  \left[  n\right]  \right]  _{q}}\right]  $. Comparison of
(\ref{ff3}) and (\ref{fm1}) leads to

\begin{theorem}
A homomorphism $\Phi^{\ast}$ of the binary group $\mathfrak{G}_{2}^{\ast}$
gives rise to a homomorphism $\Phi$ of the corresponding $n$-ary group
$\mathfrak{G}_{n}$, if $\Phi^{\ast}$ satisfies the \textquotedblleft
deformed\textquotedblright\ compatibility conditions \textrm{(\ref{ff1}%
)--(\ref{ff2})}.
\end{theorem}

The \textquotedblleft nondeformed\textquotedblright\ version ($q=1$) of this
theorem and the case of $\Phi^{\ast}$ being an isomorphism was considered in
\cite{dud/mic2}.

\bigskip

{\bf Acknowledgments}.
{The author is grateful to the Alexander von Humboldt Foundation for valuable
support
of his research stay at the University of M\"unster and to J. Cuntz for 
fruitful discussions and kind hospitality.

\end{document}